\documentclass[reqno]{amsart}
\usepackage{amsfonts}
\usepackage{amsmath}
\usepackage{amsthm}
\usepackage{amssymb}
\usepackage{latexsym}

\newif\iffinal
\finaltrue
\newtheorem*{cor}{Corollary}
\newtheorem*{lem}{Lemma}
\newtheorem*{prop}{Proposition}

\def\checkbibitem#1{\bibitem{empty}}

\theoremstyle{definition}
\newtheorem*{defn}{Definition}
\theoremstyle{definition}
\newtheorem{thm}{Theorem}
\newtheorem*{thm*}{Theorem}

\newtheorem*{rem}{Remark}

\newenvironment{pf}{\begin{proof}}{\end{proof}}
\newenvironment{pf*}{\par\normalfont\topsep6pt plus 6pt\relax\trivlist%
\item[\hskip\labelsep\itshape\proofname.]\ignorespaces}{\endtrivlist}

\newcounter{cnt}
\newenvironment{enumerit}{\begin{list}{{\hfill\rm(\roman{cnt})\hfill}}{%
\settowidth{\labelwidth}{{\rm(iv)}}\leftmargin=\labelwidth%
\advance\leftmargin by
\labelsep\rightmargin=0pt\usecounter{cnt}}}{\end{list}}

\theoremstyle{remark}

\numberwithin{equation}{section} 
\begin{document}

\newcommand{\thmref}[1]{Theorem~\ref{#1}}
\newcommand{\secref}[1]{Section~\ref{#1}}
\newcommand{\lemref}[1]{Lemma~\ref{#1}}
\newcommand{\propref}[1]{Proposition~\ref{#1}}
\newcommand{\corref}[1]{Corollary~\ref{#1}}
\newcommand{\remref}[1]{Remark~\ref{#1}}
\newcommand{\defref}[1]{Definition~\ref{#1}}
\newcommand{\er}[1]{(\ref{#1})}
\newcommand{\id}{\operatorname{id}}
\newcommand{\tensor}{\otimes}
\newcommand{\dist}{\operatorname{dist}}
\newcommand{\nc}{\newcommand}
\newcommand{\End}{\operatorname{End}}
\newcommand{\rnc}{\renewcommand}
\newcommand{\Maj}{\operatorname{Maj}}
\newcommand{\qbinom}[2]{\genfrac[]{0pt}0{#1}{#2}}
\nc{\cal}{\mathcal} \nc{\goth}{\mathfrak} \rnc{\bold}{\mathbf}
\renewcommand{\frak}{\mathfrak}
\renewcommand{\Bbb}{\mathbb}
\nc\bgamma{{\boldsymbol\gamma}}
\nc\wt{\operatorname{wt}}
\nc\Ht{\operatorname{ht}}
\newcommand{\lie}[1]{\mathfrak{#1}}
\makeatletter
\def\section{\def\@secnumfont{\mdseries}\@startsection{section}{1}%
  \z@{.7\linespacing\@plus\linespacing}{.5\linespacing}%
  {\normalfont\scshape\centering}}
\def\subsection{\def\@secnumfont{\bfseries}\@startsection{subsection}{2}%
  {\parindent}{.5\linespacing\@plus.7\linespacing}{-.5em}%
  {\normalfont\bfseries}}
\makeatother
\def\subl#1{\subsection{\iffinal\else#1\fi}\label{#1}\def\subslbl{#1}}
\def\lbl#1{\label{\subslbl.#1}}
\def\loceqref#1{\eqref{\subslbl.#1}}
\nc{\wh}[1]{\widehat{#1}}
\nc{\whbu}{\wh\bu}
\nc{\wtil}[1]{\widetilde{#1}}
\nc{\bsm}[1]{\boldsymbol{#1}}
\nc{\bpi}{{\bsm\pi}}
\nc{\bvpi}{{\bsm\varpi}}

\nc{\Cal}{\cal} \nc{\Xp}[1]{X^+(#1)} \nc{\Xm}[1]{X^-(#1)}
\nc{\on}{\operatorname} \nc{\ch}{\mbox{ch}} \nc{\Z}{{\bold Z}}
\nc{\J}{{\cal J}} \nc{\C}{{\bold C}} \nc{\Q}{{\bold Q}}
\renewcommand{\P}{{\cal P}}
\nc{\N}{{\Bbb N}} \nc\boa{\bold a} \nc\bob{\bold b} \nc\boc{\bold
c} \nc\bod{\bold d} \nc\boe{\bold e} \nc\bof{\bold f}
\nc\bog{\bold g} \nc\boh{\bold h} \nc\boi{\bold i} \nc\boj{\bold
j} \nc\bok{\bold k} \nc\bol{\bold l} \nc\bom{\bold m}
\nc\bon{\bold n} \nc\boo{\bold o} \nc\bop{\bold p} \nc\boq{\bold
q} \nc\bor{\bold r} \nc\bos{\bold s} \nc\bou{\bold u}
\nc\bov{\bold v} \nc\bow{\bold w} \nc\boz{\bold z}

\nc\ba{\bold A} \nc\bb{\bold B} \nc\bc{\bold C} \nc\bd{\bold D}
\nc\be{\bold E} \nc\bg{\bold G} \nc\bh{\bold H} \nc\bi{\bold I}
\nc\bj{\bold J} \nc\bk{\bold K} \nc\bl{\bold L} \nc\bm{\bold M}
\nc\bn{\bold N} \nc\bo{\bold O} \nc\bp{\bold P} \nc\bq{\bold Q}
\nc\br{\bold R} \nc\bs{\bold S} \nc\bt{\bold T} \nc\bu{\bold U}
\nc\bv{\bold V} \nc\bw{\bold W} \nc\bz{\bold Z} \nc\bx{\bold X}

\nc\Ca{\cal A} \nc\Cb{\cal B} \nc\Cc{\cal C} \nc\Cd{\cal D}
\nc\Ce{\cal E} \nc\Cg{\cal G} \nc\Cf{\Cal F} \nc\Ch{\cal H} \nc\Ci{\cal I}
\nc\Cj{\cal J} \nc\Ck{\cal K} \nc\Cl{\cal L} \nc\Cm{\cal M}
\nc\Cn{\cal N} \nc\Co{\cal O} \nc\Cp{\cal P} \nc\Cq{\cal Q}
\nc\Cr{\cal R} \nc\Cs{\cal S} \nc\Ct{\cal T} \nc\Cu{\cal U}
\nc\Cv{\cal V} \nc\Cw{\cal W} \nc\Cz{\cal Z} \nc\Cx{\cal X}
\nc{\hyper}[2]{\genfrac{}{}{0pt}{}{#1}{#2}}
\nc{\Ann}{\operatorname{Ann}}
\def\<#1,#2>{#1(#2)}
\nc{\indref}[1]{\ref{#1}, p.~\pageref{#1}}
\def\hatg{\wh{\lie g}}
\def\hath{\wh{\lie h}}
\def\hatbp{\wh{\mathbb P}}
\def\mbp{\mathbb P}
\newenvironment{enumspecm}[1]{\list{($#1\thecnt$)}{%
\settowidth{\labelwidth}{($#1{5}$)}%
\usecounter{cnt}\leftmargin=\labelwidth\advance\leftmargin by\labelsep%
\rightmargin=0pt\itemindent=0pt}}%
{\endlist}

\title[Filtrations and completions of positive level modules]{Filtrations and completions of certain 
positive level modules of  
 affine algebras}
\author{Vyjayanthi Chari}
\address{Department of Mathematics, University of
California, Riverside, CA 92521, U.S.A.}
\email{chari@math.ucr.edu}
\author{Jacob Greenstein}
\address{Institut de Math\'ematiques de Jussieu, Universit\'e
Pierre et Marie Curie, 175 rue du Chevaleret, Plateau 7D, F-75013
Paris, France.}
\email{greenste@math.jussieu.fr}
\date{\today}
\maketitle

This paper was motivated by an effort to understand the
representation theoretic meaning of the  results of~\cite{G,GL,NS} on 
realizations of (pseudo-)crystal bases of certain 
quantum loop modules in the framework of Littelmann's path model. 
These papers showed in particular, that one could  write  the  tensor 
product of a crystal basis of a highest weight integrable module with 
a (pseudo-)crystal basis of such a quantum loop module
as a union of highest weight crystals. The obvious and natural intepretation 
would be that the decomposition of the crystals gave rise to
a direct sum decomposition of the tensor product 
of the corresponding modules for the quantum affine algebra.  
It is however,   not very difficult to see that
 such a tensor product never
contains a copy of a highest weight module.  In addition,
the corresponding classical situation which was studied in \cite{CPnew} 
and more recently in \cite{A,Li} did not exclude the possibility that
such tensor products might in fact be irreducible.
In this paper we are able to show
that the tensor product of an integrable  highest weight representation 
with the quantum loop module associated to 
the natural representation admits a filtration
such that the  successive quotients are highest weight integrable modules  
with  multiplicity and highest weight   given by the path model.

We now  describe the main results of the paper.
In~\secref{REP} we recall some well-known properties
of highest weight modules and modules of level zero.
 We also establish several
new results on the structure of an 
irreducible finite dimensional module $V$,
and in particular introduce a function
 $n:V\to \bn$ which plays an important role in~\secref{RN}.
In~\secref{SF} we establish (\thmref{thmA})
the quantum analogue of one of the main results of~\cite{CPnew} (Theorem~4.2).
Namely, we prove that the tensor product
of a simple highest weight module with a finite 
dimensional module is simple. In this situation, we work   over the 
smaller version of the  quantum affine algebra which does not contain an 
analogue of the Euler operator. The result is the same as the classical one proved in~\cite{CPnew} but,  the absence (in general) 
of the evaluation map  and the 
non-cocommutativity of the comultiplication
in the quantum case
makes it  harder to establish. 

In the  rest of the paper we study the more complicated and interesting 
situation of the tensor product of a highest weight module~$V(\Lambda)$ with a 
quantum loop module~$L(V)$. We begin by introducing (\secref{FIL}) a filtration
$\Cv_n\supseteq\Cv_{n+1}$, $n\in\bz$ on
$V(\Lambda)\tensor L(V)$.
We prove that this   filtration is either strictly decreasing,
i.e.~$\Cv_n \supsetneq \Cv_{n+1}$ 
for all~$n\in\bz$, or trivial, i.e.~$\Cv_n=\Cv_m$ for all~$n,m\in\bz$.
We show also in this section that when $V(\Lambda)$ 
is the Verma module,
 the filtration~$\Cv_n$ is always strictly decreasing and 
that $\bigcap_{n\in\bz}\Cv_n=0$. Furthermore, for all~$n\in\bz$  the quotients 
$\Cv_n/\Cv_{n+1}$ are modules in the category $\Co$ for $\whbu_q$.
In the case when $V(\Lambda)$ is a highest weight integrable module we prove 
that $\bigcap_{n\in\bz}\Cx_n$ is zero if the filtration is strictly
 decreasing and 
that $V(\Lambda)\otimes L(V)$
 is irreducible if the filtration is trivial.

In the next two sections  we study the filtration 
$\Cx_n$ $n\in\bz$, of~$X(\Lambda)\otimes L(V)$
 where $X(\Lambda)$ is the  irreducible 
integrable module with highest weight $\Lambda$.
We  give sufficient conditions for the filtration to
 be trivial or strictly decreasing.
In the latter  case the 
 quotients $\Cx_n/\Cx_{n+1}$ are integrable
 modules in the category $\Co$ and hence isomorphic to 
finite direct sums of 
irreducible highest weight integrable modules $X(\mu)$.  These 
are related to embeddings of~$L(V)$ into
$\operatorname{Hom}_{\bc(q)}(X(\lambda),X(\mu))$ which arise in the theory of KPRV determinants
(cf.~\cite{Jo,JoT}).
 
In the last section,
we let $L(V)$ be the loop module associated to the natural representation
of the quantum affine algebra of classical type. 
and study the filtration on $X(\Lambda)\otimes L(V)$.
 This case is not covered by either of 
the sufficient conditions given in the previous sections. 
We are still  able to show that $\Cx_n\supsetneq\Cx_{n+1}$ if
 $\dim X(\Lambda)>1$.  
 We also identify the highest 
weight and multiplicites of the irreducible modules in $\Cx_n/\Cx_{n+1}$.

It follows from our results that one can complete the modules $M(\Lambda)\otimes L(V)$ (and $X(\Lambda)\otimes L(V)$ if $\Cx_n\supsetneq\Cx_{n+1}$)
with respect to the topology induced by 
the filtration. Further,  $M\Lambda)\otimes L(V)$
(resp. $X(\Lambda)\otimes L(V)$) embeds canonically into $M(\Lambda)\wh\otimes L(V)$ (resp. $X(\Lambda)\wh\otimes L(V)$). 
This can  be compared with the results of~\cite{ESV,IFR}.
 The details will apper elsewhere.

\subsection*{Acknowledgements} 
We are grateful to J.~Bernstein and O.~Schiffmann for enlightening discussions
on completions. The second author's research has been supported 
by a Marie Curie Fellowship of the European Community programme 
``Human potential'' under contract no.~HPMF-CT-2001-01132. A part
of this work was done while the second author was visiting
the Fields Institute for Research in Mathematical Sciences, Toronto,
Canada to attend  the workshop and  conference on
Infinite dimensional Lie theory and its applications. It is a 
pleasure to thank the the Institute and the  organizers of the conference 
  for their hospitality.

\section{Preliminaries}\label{P}

Throughout this paper~$\bn$ (respectively, $\bn^+$) denotes the
set of non-negative (respectively, positive) integers.

\subl{P10} Let $\frak g$ be a complex finite-dimensional simple
Lie algebra  of rank $\ell$ with a Cartan subalgebra~$\frak h$. Set
$I=\{1,2,\dots ,\ell\}$ and let $A=(d_ia_{ij})_{i,j\in I}$, where
the~$d_i$ are positive co-prime integers, be the $\ell\times \ell$ symmetrized
Cartan matrix of $\frak g$. Let $\{\alpha_i\,:\,i\in
I\}\subset\frak h^*$ (respectively $\{\varpi_i\,:\,i\in
I\}\subset\frak h^*$) be the set of simple roots (respectively  of
fundamental weights) 
of $\frak g$ with respect to $\frak h$. Let~$\theta$
be the highest root of~$\lie g$. 
As
usual, $Q$ (respectively, $P$) denotes the root (respectively,
weight) lattice of $\frak g$. Let~$P^+=\sum_{i\in I}\bn\varpi_i$
be the set of dominant weights and set~$Q^+=\sum_{i\in
I}\bn\alpha_i$. Given~$\gamma=\sum_{i\in I} k_i\alpha_i\in Q^+$, set~$\Ht\gamma=\sum_{i\in I}k_i$.
Let $W$ be the Weyl group of $\lie g$ and let $s_\alpha\in W$ denote 
the relection with respect to the root $\alpha$.  It is well-known that~$\lie h^*$  admits a
non-degenerate symmetric $W$-invariant  bilinear form which will be denoted
by~$(\cdot\,\vert\,\cdot)$. We assume that~$(\alpha_i\,\vert\,
\alpha_i)=d_i a_{ij}$ for all~$i,j\in I$. Given a root~$\beta$ of~$\lie g$,
denote by~$\beta^\vee\in\lie h$ the corresponding co-root. 

\subl{P15}
Let
$$\widehat{\frak g}=\frak g\otimes\bc[t,t^{-1}]\oplus\bc
c\oplus \bc d
$$
be the untwisted extended affine algebra
associated with~$\frak g$ and let $\wh{A}=(d_ia_{ij})_{i,j\in\wh
I}$, where~$\wh I=I\cup\{0\}$ be the extended symmetrized Cartan
matrix. Set $\widehat{\frak h}=\frak h\oplus \bc c\oplus \bc d$.

From now on we identify $\lie h^*$ with the subspace of $\wh{\lie
h}^*$ consisting of elements which are zero on $c$ and $d$.
Define~$\delta\in \wh{\frak h}^*$ by $$\delta(\frak h\oplus \bc
c)=0,\qquad\delta(d)=1.$$ 
Set $\alpha_0=\delta-\theta$. Then
$\{\alpha_i\,:\,i\in\wh I\}$ is a set of simple roots for
$\wh{\lie g}$ with respect to~$\wh{\lie h}$, $\alpha_0^\vee=
c-\theta^\vee$ and $\delta$ generates
its imaginary roots. 

Let~$\wh W$ be the Weyl group of~$\hatg$. 
The bilinear form on $\frak h^*$ extends to a
$\wh W$-invariant  bilinear form on $\hath^*$ which we
continue to denote by~$(\cdot\,\vert\,\cdot)$. One has $(\delta\,
\vert\,\alpha_i)=0$ and $(\alpha_i\,\vert\,\alpha_j)=d_i a_{ij}$,
for all~$i,j\in\wh I$. Define a set of fundamental
weights~$\{\omega_i\,:\, i\in \wh I\}\subset\widehat{\frak h}^*$
of~$\widehat{\frak g}$ by the conditions~$(\omega_i|\alpha_j)=d_i
\delta_{i,j}$ and $\omega_i(d)=0$ for all~$i,j\in\wh I$. Let~$\widehat
P= \sum_{i\in\wh I}\bz\omega_i\oplus\Z\delta$ (respectively,
$\widehat P^+=\sum_{i\in\wh I}\bn\omega_i\oplus\Z\delta$) be the
corresponding set of integral (respectively, dominant) weights.
We have $\varpi_i=\omega_i-a_i^\vee
\omega_0$ where~$a_i^\vee$ is the coefficient of~$\alpha_i^\vee$ 
in~$\theta^\vee$. Identify~$P$ with the free abelian subgroup
of~$\wh P$ generated by the~$\varpi_i$, $i\in I$.
Denote by~$\wh Q$ 
the root lattice of~$\hatg$ and set~$\wh Q^+=\sum_{i\in\wh I} \bn\alpha_i$. 
Given~$\wh\gamma=\sum_{i\in \wh I}k_i\alpha_i\in\wh Q^+$, set~$\Ht\wh\gamma=\sum_{i\in\wh I} k_i$.
Given $\lambda,\mu\in
\widehat{P}^+$ (respectively, $\lambda,\mu\in P^+$) we say that
$\lambda\le \mu$ if $\mu-\lambda\in \widehat{Q}^+$ (respectively,
$\mu-\lambda\in Q^+$). For all~$\lambda\in\wh P$ set~$\lambda_i=\lambda(\alpha_i^\vee)$, $i\in\wh I$.

\subl{P20} Let $q$ be
an indeterminate and let $\bc(q)$ be the field of rational
functions in $q$ with complex coefficients.  For $r,m\in\bn$,
$m\ge r$, define
\begin{equation*}
[m]_q=\frac{q^m -q^{-m}}{q -q^{-1}},\qquad [m]_q!
=[m]_q[m-1]_q\ldots [2]_q[1]_q,\qquad \qbinom{m}{r}_q =
\frac{[m]_q!}{[r]_q![m-r]_q!}.
\end{equation*} 
For $i\in\wh I$, set $q_i=q^{d_i}$ and
$[m]_i=[m]_{q^i}$. 

The quantum affine algebra~$\widehat\bu_q(\frak
g)$ (cf.~\cite{Be,BCP,Dr,J}) associated to $\frak g$, which will
be further denoted as~$\widehat\bu_q$, is an associative algebra
over $\bc(q)$ with generators $x_{i,r}^{{}\pm{}}$, $h_{i,k}$,
$K_i^{ {}\pm 1}$, $C^{\pm 1/2}$, $D^{\pm
1}$, where~$i\in I$, $k,r\in \bz$, $k\not=0$, and the
following defining relations
\begin{gather*}
  \text{$C^{{}\pm{1/2}}$ are central,}\\
K_iK_i^{-1} = K_i^{-1}K_i
  =1,\quad C^{1/2}C^{-1/2} =C^{-1/2}C^{1/2} =1, \qquad DD^{-1}=D^{-1}D=1\\
K_iK_j =K_jK_i,\quad D K_i = K_i D\\
\quad K_ih_{j,r} =h_{j,r}K_i,\quad Dh_{j,r}D^{-1}=q^r
h_{j,r},
\\
K_ix_{j,r}^\pm K_i^{-1} = q_i^{{}\pm a_{ij}}x_{j,r}^{{}\pm{}},\qquad
D x_{j,r}^{\pm} D^{-1} = q^r x_{j,r}^{\pm},
\\
  [h_{i,r},h_{j,s}]=\delta_{r,-s}\,\frac1{r}[ra_{ij}]_i\,\frac{C^r-C^{-r}}
  {q_j^{}-q_j^{-1}},
\\
[h_{i,r} , x_{j,s}^{{}\pm{}}] =
  \pm\frac1r[ra_{ij}]_i\, C^{{}\mp {|r|/2}}x_{j,r+s}^{{}\pm{}},
\\
\intertext{}
  x_{i,r+1}^{{}\pm{}}x_{j,s}^{{}\pm{}} -q_i^{{}\pm
    a_{ij}}x_{j,s}^{{}\pm{}}x_{i,r+1}^{{}\pm{}}=q_i^{{}\pm
    a_{ij}}x_{i,r}^{{}\pm{}}x_{j,s+1}^{{}\pm{}}
  -x_{j,s+1}^{{}\pm{}}x_{i,r}^{{}\pm{}},
\\
[x_{i,r}^+ ,
  x_{j,s}^-]=\delta_{i,j}\,\frac{ C^{(r-s)/2}\,\psi_{i,r+s}^+ -
    C^{-(r-s)/2}\, \psi_{i,r+s}^-}{q_i^{} - q_i^{-1}},
\\
\sum_{\pi\in\Sigma_m}\sum_{k=0}^m(-1)^k\qbinom{m}{k}_{i}
  x_{i, r_{\pi(1)}}^{{}\pm{}}\cdots x_{i,r_{\pi(k)}}^{{}\pm{}}
  x_{j,s}^{{}\pm{}} x_{i, r_{\pi(k+1)}}^{{}\pm{}}\cdots
  x_{i,r_{\pi(m)}}^{{}\pm{}} =0,\qquad\text{if $i\ne j$},
\end{gather*}
for all sequences of integers $r_1,\ldots, r_m$, where $m
=1-a_{ij}$, $\Sigma_m$ is the symmetric group on $m$ letters, and
the $\psi_{i,r}^{{}\pm{}}$ are determined by equating powers of
$u$ in the formal power series
$$
\sum_{r=0}^{\infty}\psi_{i,\pm
r}^{{}\pm{}}u^{{}\pm r} = K_i^{{}\pm 1}
\exp\Big(\pm(q_i^{}-q_i^{-1})\sum_{s=1}^{\infty}h_{i,\pm s}
u^{{}\pm s} \Big).
$$
The
subalgebra of $\whbu_q$ generated by the elements
$x_{i,0}^\pm$, $i\in I$ is isomorphic to the quantized enveloping
algebra $U_q(\lie g)$ of $\lie g$.

Let~$\whbu_q(\gg)$ (respectively, $\whbu_q(\ll)$) 
be the subalgebra of~$\whbu_q$
generated by the~$x_{i,s}^+$ (respectively, by the $x_{i,s}^-$) for all~$i\in I$, $s\in \bz$.
Given~$r\in\bz$, let~$\whbu_q^r(\gg)$ (respectively, $\whbu_q^r(\ll)$) be the
subalgebra of~$\whbu_q(\gg)$ (respectively, of $\whbu_q(\ll)$) generated by
the~$x_{i,s}^+$ (respectively, by the~$x_{i,s}^-$) 
for all~$i\in I$ and for all~$s\ge r$.
Furthermore, let~$\whbu_q(0)$ (respectively, $\whbu_q^r(0)$) be the subalgebra of~$\whbu_q$
generated by the~$h_{i,s}$, for all~$i\in I$ and for all~$s\in\bz$ (respectively,
for all~$s\ge r$), $s\not=0$. Finally, let $\whbu_q^\circ$ be the subalgebra generated by the $K_i^{\pm 1}$, 
$i\in I$, $D^{\pm1}$ and~$C^{\pm1/2}$. 

\subl{P23}
Define a~$\bz$-grading on~$\whbu_q$ by setting~$\deg x_{i,r}^{\pm}=r$,
$\deg h_{i,k}=k$ for all~$j\in I$ and for all~$r\in\bz$, $k\in\bz\setminus
\{0\}$ and~$\deg K_i = \deg D=\deg C^{\pm1/2}=0$ for all~$i\in I$.
Equivalently, we say that~$x\in\whbu_q$ is homogeneous of degree~$k=\deg x$
if~$DxD^{-1}=q^k x$. Given~$z\in\bc(q)^\times$, 
let~$\phi_z$ be the automorphism of~$\whbu_q$ defined by extending~$\phi_z(x)=
z^{\deg x} x$ for~$x\in\whbu_q$ homogeneous.

On the other hand, the algebra~$\whbu_q$ is graded by the root 
lattice~$\wh Q$, the elements~$x_{i,r}^\pm$, $i\in I$, $r\in\bz$ being of weight~$r\delta
\pm\alpha_i$, the~$h_{i,k}$, $i\in I$, $k\in\bz\setminus\{0\}$ being of weight~$k\delta$ and the other generators
being of weight zero. Given~$\nu\in\wh Q$, we denote the corresponding
weight subspace of~$\whbu_q$ by~$(\whbu_q)_\nu$. Observe also that
if~$x\in(\wh\bu_q)_{r\delta+\gamma}$, $\gamma\in Q$, $r\in\bz$ then~$\deg x=r$.

\subl{P25} We will also need another presentation
of~$\whbu_q$. Namely, after~\cite{Be,J}, the algebra
$\whbu_q$ is isomorphic to an associative $\bc(q)$-algebra
generated by $E_i, F_i, K_i^{{}\pm 1}\,:\,i\in\wh I$, $D^{\pm 1}$
and central elements~$C^{\pm1/2}$ satisfying the following
relations:
\begin{gather*}
\text{$C=K_0\displaystyle\prod_{i\in I} K_i^{a_i}$, where~$\theta=
\displaystyle\sum_{i\in I} a_i\alpha_i$, $a_i\in\N^+$}
\\
\intertext{}
K_iE_j K_i^{-1}=q_i^{ a_{ij}}E_j,\qquad K_iF_j
K_i^{-1}=q_i^{-a_{ij}}F_j,\\ DE_j D^{-1}=q^{
\delta_{j0}}E_j,\qquad DF_j D^{-1}=q^{-\delta_{j0}}F_j,\\ 
[E_i,F_j]=\delta_{i,j}\,\frac{K_i-K_i^{-1}}{q_i^{}-q_i^{-1}},\\
  \sum_{r=0}^{1-a_{ij}}(-1)^r\qbinom{1-a_{ij}}{r}_i
(E_i)^rE_{j}(E_{i})^{1-a_{ij}-r}=0\
  \qquad \text{if $i\ne j$},\\
\sum_{r=0}^{1-a_{ij}}(-1)^r\qbinom{1-a_{ij}}{r}_i
(F_{i})^rF_{j}(F_{i})^{1-a_{ij}-r}=0\
  \qquad \text{if $i\ne j$}.
\end{gather*}
The element ~$E_i$ (respectively, $F_i$), $i\in I$ corresponds
to~$x_{i,0}^+$ (respectively, $x_{i,0}^-$). In particular, the
elements $E_i,F_i,K_i^{\pm1}\,:\,i\in I$ generate a subalgebra of~$\whbu_q$
isomorphic to~$U_q(\frak g)$.

Let $\wh{\bu}_q^+$ (respectively. $\whbu_q^-$) be the
$\bc(q)$-subalgebra of $\whbu_q$ generated by the~$E_i$ (respectively, by the $F_i$), 
$i\in \wh{I}$. Let~$\whbu_q'$ be the  subalgebra of~$\whbu_q$
generated by the~$E_i, F_i, K_i^{\pm1}$, $i\in\wh I$ and
by~$C^{\pm1/2}$.

We will need the following
result which was established in~\cite{BCP}.
\begin{prop}\label{triangle} We have
$\whbu_q^+\subset \whbu_q^0(\ll)\whbu_q^0(0)\whbu_q^0(\gg)$ and~$\whbu_q^r(\ll),\whbu_q^s(\gg)\subset \whbu_q^+$
for all~$r\in\bn^+$, $s\in\bn$.
\end{prop}

\subl{P30}
It is well-known that $\whbu_q$ is a Hopf algebra over
$\bc(q)$ with  the co-multi\-pli\-cation being given in terms of
generators $E_i,F_i,K_i^{\pm1} \,:\,i\in\wh I$ by the following
formulae $$ \Delta(E_i)=E_i\tensor 1+K_i\tensor E_i,\qquad
\Delta(F_i)=F_i\tensor K_i^{-1} + 1\tensor F_i, $$
the~$K_i^{\pm1}$, $D^{\pm1}$, $C^{\pm1/2}$ being group-like.
Notice that $\whbu_q'$ is a Hopf subalgebra of $\whbu_q$.
Let~$\wh\Cu_q^+$ (respectively, $\wh\Cu_q^-$) be the subalgebra of~$\wh\bu_q$ and~$\wh\bu_q'$ generated
by the~$E_i$ (respectively, by the~$F_i$) and by the~$K_i^{\pm1}$, $i\in\wh I$. 
Obviously, the $\wh\Cu_q^\pm$ are Hopf algebras and~$\wh\bu_q^\pm\subset\wh\Cu_q^\pm$.

Although explicit formulae for the co-multiplication on generators
~$x_{i,r}^\pm$, $h_{i,r}$ are not known, we have the following
partial results \cite{Da} which are enough for this paper.
\begin{lem}\label{comul}
For $i\in I$,  $r\in\bn$, $s\in\bn^+$, we have
\begin{alignat}{3}
&\Delta(h_{i,s})=h_{i,s}\otimes 1+1\otimes h_{i,s} &+\,&\text{\em terms
in~$\whbu_q^\circ((\whbu_q^+)_+\otimes (\whbu_q^+)_+)$},\lbl{10a}\\
&\Delta(x^+_{i,r})=x^+_{i,r}\otimes 1+ K_i\otimes x^+_{i,r} &+
\,&\text{\em terms in~$\whbu_q^\circ((\whbu_q^+)_+\otimes(\whbu_q^0(\gg))_+)$},
\lbl{10b}\\
&\Delta(x^-_{i,s})=x^-_{i,s}\otimes K_i+ 1\otimes x^-_{i,s} &+
\,&\text{\em terms in~$\whbu_q^\circ((\whbu_q^+)_+\otimes (\whbu_q^1(\ll))_+)$},
\lbl{10c}
\end{alignat}
where~$(\whbu_q^+)_+$ denotes the augmentation ideal of~$\whbu_q^+$.
\end{lem}
%
For $i\in I$,  set
$$
P^\pm_i(u)=\exp\Big(-\sum_{k=1}^\infty
\frac{q^{\pm k}h_{i,\pm k}}{[k]_i}u^k\Big).
$$ 
Let $P_{i,\pm r}$
be the coefficient of $u^r$ in $P_i^\pm (u)$. It is easy to see that the
elements $h_{i,r}$ belong to the subalgebra of $\whbu_q$ generated by the
elements $P_{i,r}$, $i\in I$, $r\in\bz$. Further, 
one can deduce from~\lemref{P30} as in~\cite{BCP} that, for all~$s\in\bn$,
\begin{equation}\lbl{10d}
\Delta(P_{i,s})=\sum_{r=0}^s P_{i,s-r}\otimes P_{i,r}+\, \text{terms
in~$\whbu_q^\circ((\whbu_q^+)_+)\otimes (\whbu_q^+)_+$})
\end{equation}

\section{The modules $M(\Lambda)$, $X(\Lambda)$, $V(\bpi)$ and $L(V(\bpi))$.}
\label{REP}

In this section we recall the definition and some properties of
several families of integrable modules for $\whbu_q$ and
$\whbu_q^\prime$. For modules of level zero we also establish
some results which we need in later sections.

\subl{REP0}
A $\whbu_q$-module $M$ is said to be of type~$1$ 
if~$M=\bigoplus_{\mu\in \wh P}M_\mu$,
where
$$
M_\mu=\{m\in M\,:\,K_im=q_i^{\mu(\alpha_i^\vee)}m,\,
\forall\ i\in\wh I,\, Dm=q^{\mu(d)} m\}.
$$
 Type~$1$-modules for~$\whbu_q'$ are defined in the obvious way.
 If~$m\in M_\mu\setminus\{0\}$, 
we say that~$m$ is of weight~$\mu$ and write~$\wt m=\mu$.
Set~$\Omega(M)=\{\nu\in\wh P\,:\, M_\nu\not=0\}$.

A $\whbu_q$- or a~$\whbu_q'$-module $M$~of type~$1$
is said to be integrable if 
the elements $E_i$, $F_i$, $i\in\wh I$ act locally
nilpotently on $M$. 
Evidently, a $\whbu_q$-module~$M$ can be viewed as a 
$\whbu_q'$-module~$M'$ and~$M'_\nu=\bigoplus_{r\in\bz} M_{\nu+r\delta}$.

\subl{REP10} Let $\cal{O}$ be the category of $\whbu_q$-modules
satisfying the following properties. A $\whbu_q$-module~$M$ 
is an object in~$\Co$ if and only if
\begin{enumerit}
\item $M$ is a module of type~$1$ and~$\dim M_\mu<\infty$ for all~$\mu\in
\wh P$.
\item The set~$\Omega(M)$ is contained in the set~$\bigcup_{k=1}^r
\{\lambda_k-\wh\gamma\,:\, \wh\gamma\in\wh Q^+\}$ for some~$r\in\bn^+$
and for some~$\lambda_k\in\wh P$.
\end{enumerit}

Given $\Lambda\in\wh P$, let $M(\Lambda)$ denote
the Verma module of highest weight~$\Lambda$. It is generated as a $\whbu_q$-module 
by an element $m_\Lambda$ of weight~$\Lambda$ with
defining relation
$$
(\whbu_q^+)_+ m_\Lambda=0
$$
It is well-known that $M(\Lambda)$ has a unique simple quotient which we denote by $X(\Lambda)$.
Let~$v_\Lambda$ be the canonical image of~$m_\Lambda$ in~$X(\Lambda)$.

The next result is well-known and follows immediately from~\cite{Kac} and~\cite{Lus}.
\begin{prop}
\begin{enumerit}
\item For all $\Lambda\in\wh P$, $M(\Lambda)\in\cal{O}$ and
is a free $\whbu_q^-$-module. In 
particular $\Ann_{\whbu_q^-} m_\Lambda=0$ and  $\Omega(M(\Lambda))\subset \Lambda-\wh Q^+$.
\item  For all $\Lambda\in\wh P^+$,  $X(\Lambda)$
is an   integrable $\whbu_q$-module in the category
$\cal{O}$ and is generated as a $\whbu_q$-module by the element
$v_\Lambda$.  
Moreover, $\Ann_{\whbu_q^-} v_\Lambda=\sum_{i\in\wh I}
\whbu_q^- F_i^{\Lambda(\alpha_i^\vee)+1}$.
\item 
Let $M\in\cal{O}$ be integrable. Then $M$ is
is isomorphic to a finite direct sum of modules of the form
$X(\Lambda)$, $\Lambda\in\wh P^+$. In particular, if~$M\in\Co$
is simple and integrable, 
then it is isomorphic to~$X(\Lambda)$ for some~$\Lambda\in\wh P^+$
\end{enumerit}
\end{prop}
Regarded as modules for $\whbu_q^\prime$, the $X(\Lambda)$ remain
simple, although they no longer have finite-dimensional
weight spaces. Indeed, by~\cite{Kac} if~$\lambda\in \Omega(X(\Lambda))$
then~$\lambda-n\delta\in \Omega(X(\Lambda))$ for all~$n\in\bn$
and we can  write
\begin{equation}{\lbl{20}}
X(\Lambda)=\bigoplus_{\gamma\in Q^+}\bigoplus_{n\in\bn} X(\Lambda)_{\Lambda-\gamma-n\delta}.
\end{equation} 
Obviously, for~$\gamma\in Q^+$ fixed, $\bigoplus_{n\in\bn} X(\Lambda)_{\Lambda-\gamma-n\delta}$
is a weight space of~$X(\Lambda)$ viewed as a $\whbu_q'$-module.
Observe also that $$M(\Lambda)\cong M(\Lambda+r\delta),\qquad
X(\Lambda)\cong X(\Lambda+r\delta)$$
as~$\whbu_q'$-modules for all~$r\in\bz$. 

\subl{REP18} 
The next important family of modules we consider is that of
the irreducible finite-dimensional
representations $V(\bpi)$ of $\whbu_q^\prime$. Let
$\bpi=(\pi_i(u))_{i\in I}$ be an $\ell$-tuple of polynomials~ with
coefficients in $\bc(q)$ and with constant term 1.
Set~$\lambda_\bpi:=\sum_{i\in I}
(\deg \pi_i)\varpi_i\in P$.
Let $W(\bpi)$ be the
$\whbu_q^\prime$-module generated by an element $v_\bpi$
satisfying
\begin{gather*}
x^+_{i,r}v_\bpi =0,\quad
(x^-_{i,r})^{\lambda_\bpi(\alpha_i^\vee)+1}v_\bpi=0,\\
K_i v_\bpi=q_i^{\lambda_\bpi(\alpha_i^\vee)}v_\bpi,\quad C v_\bpi=v_\bpi,\\
P_{i,\pm s}v_\bpi=\pi^\pm_{i,s}v_\bpi,
\end{gather*} 
for all $i\in I$, $r\in\bz$,$s\in\bn$ where~$\pi_i^\pm(u)=\sum_{s} \pi_{i,s}^\pm
u^s$ and
$$
\pi_i^+(u)=\pi_i(u), \qquad
\pi_i^-(u)=
u^{\deg\pi_i} \frac{\pi_i(u^{-1})}{\pi^+_{i,\deg\pi_i}}.
$$ 
The following proposition was proved in~\cite{CPweyl}.
\begin{prop} 
\begin{enumerit} 
\item The $\whbu_q^\prime$-modules $W(\bpi)$ are finite-dimensional.
\item $W(\bpi)=\whbu_q(\ll)v_\bpi$. In particular, $\Omega(W(\bpi))\subset
\lambda_\bpi-Q^+$.
\item 
$\dim W(\bpi)_{\lambda_\bpi}=\dim W(\bpi)_{w_\circ\lambda_\bpi}=1$,
where~$w_\circ$ is the longest element of the Weyl group of~$\lie
g$. Let $v_\bpi^*$ be a non-zero element in
$W(\bpi)_{w_\circ\lambda_\bpi}$. Then 
$$
x^-_{i,r}v_\bpi^*=0,\quad (x^+_{i,r})^{-(w_\circ\lambda_\bpi)(\alpha_i^\vee)
+1}v_{\bpi}^*=0,\quad
W(\bpi)=\whbu(\gg)v_\bpi^*,$$ and $\Omega(W(\bpi))\subset
w_\circ\lambda_\bpi+Q^+.$
\item 
$W(\bpi)$ has a unique simple quotient $V(\bpi)$ and all simple
finite dimensional $\whbu_q'$-modules are obtained that way.
\item 
Denote the images of the elements  $v_\bpi$, $v_\bpi^*$ in
$V(\bpi)$ by the same symbols. Then
\begin{gather*}
V(\bpi)=\whbu_q(\ll)v_\bpi,\quad
x^+_{i,r}v_\bpi =0,\quad (x_{i,r}^-)^{\lambda_{\bpi}(\alpha_i^\vee)+1}
v_{\bpi}=0,\\ 
\Omega(V(\bpi))\subset \lambda_\bpi-Q^+,
\end{gather*}
and analogous statements hold for $v_\bpi^*$.
\end{enumerit}
\end{prop}
Given~$z\in\bc^\times$ and an~$\ell$-tuple of polynomials~$\bpi=(\pi_i(u))_{i\in I}$, 
one can introduce on~$V(\bpi)$ another $\whbu_q'$-module structure
by twisting the action by automorphism~$\phi_z$. 
Then~$\phi_z^*V(\bpi)\cong V(\bpi_z)$ where~$\bpi_z=(\pi_i(zu))_{i\in I}$.

\subl{REP20}
We now establish some facts about~$V(\bpi)$ which will be
needed later.
\begin{lem} Let $\bpi$ be an $\ell$-tuple of polynomials with
coefficients in $\bc(q)$ and constant term 1. Let~$k_i$
(respectively, $k^*_j$), $i,j\in I$, be the dimension
of~$V(\bpi)_{\lambda_\bpi-\alpha_i}$ (respectively,
of~$V(\bpi)_{w_\circ\lambda_\bpi+\alpha_j}$). Suppose
that~$k_i,k_j^*>0$. Then
\begin{enumerit}
\item $\{x_{i,s}^- v_{\bpi},\dots,x_{i,s+k_i-1}^-
v_{\bpi}\}$ is a basis of~$V(\bpi)_{\lambda_\bpi-\alpha_i}$
for all~$s\in\bz$.
\item $\{x_{j,s}^+ v_{\bpi}^*,
\dots,x_{j,s+k^*_j-1}^+v_{\bpi}^*\}$ is a basis of~$V(\bpi)_{w_\circ
\lambda_\bpi+\alpha_j}$ for all~${s\in\bz}$.
\end{enumerit}
\end{lem}
\begin{pf}
We prove only~(i), the proof of~(ii) being similar. Since~$V(\bpi)=
\whbu_q(\ll)v_\bpi$, the elements~$x_{i,k}^- v_{\bpi}$, $k\in\bz$
span~$V(\bpi)_{\lambda_{\bpi}-\alpha_i}$. Next, observe
that~$x_{i,k}^- v_{\bpi}\not=0$ for all~$k\in\bz$. Indeed,
if~$x_{i,n}^- v_{\bpi}=0$ for some~$n\in\bz$ then,
since~$v_{\bpi}$ is an eigenvector for the~$h_{i,s}$, $s\in\bz$ we
get, using the defining relations of~$\whbu_q$ 
$$ 
0=h_{i,s}x_{i,n}^- v_{\bpi}=-\frac{[2s]_i}{s}\, x_{i,n+s}^- v_{\bpi}. 
$$ 
It follows that~$x_{i,k}^- v_{\bpi}=0$ for all~$k\in\bz$. Therefore,
$V(\bpi)_{\lambda_\bpi-\alpha_i}=0$, which is a contradiction.

It remains to prove that
the set~$\{x_{i,s}^- v_{\bpi},\dots,x_{i,s+k_i-1}^-v_{\bpi}\}$
is linearly independent for all~$s\in\bz$. If~$k_i=1$ then,
since~$x_{i,s}^- v_{\bpi}\not=0$ for all~$s\in\bz$, there is nothing
to prove. Assume that~$k_i>1$ and that
$\sum_{r=0}^{k_i-1} a_r x_{i,r+s}^- v_{\bpi}=0$ for some~$a_r\in\bc(q)$,
$r=1,\dots,k$ and for some~$s\in\bz$.
Applying~$h_{i,m}$ as above we conclude
that
\begin{equation}\lbl{10}
\sum_{r=0}^{k_i-1} a_r x_{i,r+s}^- v_{\bpi}=0,\qquad \forall s\in\bz.
\end{equation}
Let~$r_1$ (respectively, $r_2$)
be the
minimal (respectively, the maximal) $r$, $0\le r \le k_i-1$
such that~$a_r\not=0$. Then, using~\loceqref{10} with~$s=0$ and~$s=1$, we obtain
$$
x_{i,r_1}^- v_{\bpi} = - a_{r_1}^{-1} \sum_{r=r_1+1}^{r_2} a_r x_{i,r}^- v_{\bpi},
\qquad x_{i,r_2+1}^- v_{\bpi} = -a_{r_2}^{-1} \sum_{r=r_1+1}^{r_2} a_{r-1} x_{i,r}^- v_{\bpi}.
$$
Observe that both sums contain at least one non-zero term. Then
it follows by induction on~$r_1-k$ (respectively, on~$k-r_2$) from the above formulae and~\loceqref{10}
that the~$x_{i,k}^- v_{\bpi}$ lie in the linear span of vectors $x_{i,r_1+1}^- v_{\bpi},\dots,
x_{i,r_2}^- v_\bpi$ for all~$k<r_1$ (respectively, for all~$k>r_2$).
Therefore, $\dim V(\bpi)_{\lambda_{\bpi}-\alpha_i}<k_i$
which is a contradiction.
\end{pf}

\subl{REP25}
\begin{lem}
Define
\begin{align*}
k(\bpi)&=\min_{i\in I} \{\dim V(\bpi)_{\lambda_{\bpi}-\alpha_i}\,:\,
V(\bpi)_{\lambda_{\bpi}-\alpha_i}\not=0\}.\\
\intertext{Then}
k(\bpi)&=\min_{i\in I} \{\dim V(\bpi)_{w_\circ\lambda_{\bpi}+\alpha_i}\,:\,
V(\bpi)_{w_\circ\lambda_{\bpi}+\alpha_i}\not=0\}.
\end{align*}
\end{lem}
\begin{pf}
Let~$k^*=\min_{i\in I} \{\dim V(\bpi)_{w_\circ\lambda_{\bpi}+\alpha_i}\,:\,
V(\bpi)_{w_\circ\lambda_{\bpi}+\alpha_i}\not=0\}$. Choose~$i\in I$
such that 
$k(\bpi)=\dim V(\bpi)_{\lambda_\bpi-\alpha_i}$ for some~$i\in I$.
Since~$V(\bpi)$ is an integrable $\wh\bu_q'$-module, 
its character is~$W$-invariant (cf. say~\cite{L}). 
Since~$w_\circ \alpha_i = -\alpha_j$ for some~$j\in I$ we conclude that
$\dim V(\bpi)_{w_\circ\lambda_\bpi+\alpha_j}=k(\bpi)$ and so~$k^*\le k(\bpi)$.
A similar argument shows that~$k(\bpi)\le k^*$.
\end{pf}

\subl{REP28} 
\begin{lem}
For any~$v\in V(\bpi)$, $E_0^{\lambda_{\bpi}(\theta^\vee)+1} v=0=
F_0^{\lambda_{\bpi}(\theta^\vee)+1}v$.
\end{lem}
\begin{pf}
Since~$V(\bpi)$ is finite dimensional, it decomposes, uniquely, as
a direct sum of simple finite dimensional highest weight
modules~$V(\lambda)$ over~$U_q(\lie g)$ with~$\lambda\in
\lambda_{\bpi}-Q^+$. Therefore, in order to prove the assertion it
is sufficient to show that, if~$\lambda- \gamma$, $\gamma\in Q^+$
is a weight of~$V(\lambda)$,
then~$\mu=\lambda-\gamma-(\lambda(\theta^\vee)+1)\theta$ is not a
weight of~$V(\lambda)$. Indeed, otherwise, since the formal
character of~$V(\lambda)$ is $W$-invariant, $s_\theta\mu=
\lambda-\gamma+{(\gamma(\theta^\vee)+1)}\theta$ is also a weight
of~$V(\lambda)$. It follows that~$\gamma'=\gamma-
(\gamma(\theta^\vee)+1)\theta\in Q^+$.

Let~$J=\{i\in I\,:\, \alpha_i(\theta^\vee)>0\}$ and observe
that~$J$ is not empty. Write~$\gamma=\sum_{i\in I} n_i \alpha_i$.
Suppose first that~$\gamma(\theta^\vee)=0$. Then~$n_i=0$ for
all~$i\in J$. Yet~$\theta=\sum_{i\in I} a_i \alpha_i$ and~$a_i>0$
for all~$i\in I$. It follows that the~$\alpha_i$, $i\in J$ occur
in~$\gamma'=\gamma-\theta$ with strictly negative coefficients.
Therefore, $\gamma'\notin Q^+$ which is a contradiction.

Finally, suppose that~$\gamma(\theta^\vee)>0$. Then there
exists~$i\in J$ such that~$n_i\not=0$. It follows that~$\alpha_i$
occurs in~$ (\gamma(\theta^\vee)+1)\theta$ with the coefficient at 
least~$a_i(n_i+1)>n_i$. Thus, $\alpha_i$ occurs in~$\gamma'$ with
a negative coefficient and so~$\gamma'\notin Q^+$.

A similar argument shows that~$F_0^{-w_\circ\lambda_\bpi(\theta^\vee)+1}
v=0$ for all~$v\in V(\bpi)$. It remains to observe 
that~$-w_\circ\lambda_\bpi(\theta^\vee)=\lambda_{\bpi}(\theta^\vee)$.
\end{pf}

\subl{REP30}
\begin{prop}
For all~$r\in\bn$,
$V(\bpi)=\whbu_q^r(\ll)v_{\bpi}=
\whbu_q^r(\gg) v_{\bpi}^*$. 
\end{prop}
\begin{pf}
It is sufficient to prove the statement for~$v_{\bpi}$, the proof
of the other one being similar. Recall that all weights of~$V(\bpi)$ are
of the form~$\lambda_{\bpi}-\gamma$, $\gamma\in Q^+$.
We prove by induction on~$\Ht\gamma$
that~$$V(\bpi)_{\lambda_\bpi-\gamma}\subset\
\whbu_q^r(\ll)v_{\bpi},\qquad \forall\,r\in\bn.$$

If~$\Ht\gamma=0$ then there is nothing to prove. Assume
that~$\Ht\gamma=1$, that is~$\gamma=\alpha_i$ for some~$i\in I$.
Then~$k=\dim V(\bpi)_{\lambda_\bpi- \alpha_i}>0$ and
by~\lemref{REP20}, $V(\bpi)_{\lambda_\bpi-\alpha_i}$ is spanned
by~$x_{i,r}^- v_{\bpi},\dots, x_{i,r+k-1}^- v_{\bpi}$ for
all~$r\ge 0$. In particular, $$V(\bpi)_{\lambda_\bpi-\alpha_i}
\subset \whbu_q^r{(\ll)}v_\bpi$$ for all~$r\ge 0$. 

Suppose that~$v\in V(\bpi)_{\lambda_\bpi-\gamma}$ with~$\Ht\gamma \ge 1$
and that $v\in \whbu_q^r(\ll)v_{\bpi}$ for all~$r\in\bn$. Fix some~$r\in\bn$.
For the inductive
step, it suffices  to prove that~$x_{i,k}^- v\in\whbu_q^r(\ll)v_{\bpi}$
for all~$k\in\bz$ and for all~$i\in I$.
We may assume, without loss of generality,
that~$v=x_{j,n}^- w$ for some~$w\in\whbu_q^{r+1}(\ll) v_{\bpi}$, $j\in I$
and~$n>r$.
It follows from the defining relations of the algebra~$\whbu_q$ that
$$
x_{i,k}^- x_{j,n}^- w = q_i^{a_{ij}} x_{j,n}^- x_{i,k}^- w+
q_i^{a_{ij}} x_{i,k+1}^- x_{j,n-1}^- w - x_{j,n-1}^- x_{i,k+1}^-w.
$$
If~$k\ge r-1$ then all terms in the right hand side lie
in~$\whbu_q^r(\ll)v_{\bpi}$ by the assumption on~$w$ and by the
induction hypothesis. Then it follows from
the above formula by induction on~$r-k$ that~$x_{i,k}^- x_{j,n}^- w
\in\whbu_q^r{(\ll)}v_{\bpi}$ for all~$k<r$ and the proposition is proved.\end{pf}
\begin{cor} 
We have $V(\bpi)=\whbu_q^+v_\bpi=\whbu_q^+ v_\bpi^*$. 
\end{cor}
\begin{pf}
This follows immediately from the above and \propref{P25}.
\end{pf}

\subl{REP35} 
As a consequence of~\propref{REP30}, we can define a map~$n:V(\bpi)\to \bn$
in the following way. Given~$v\in V(\bpi)$, $v\not=0$ let~$n(v)$ be the minimal~$r\in\bn$ such
that~$v$ can be written as a linear combination of homogeneous elements of~$\wh\bu_q'$ of
degree~$\le r$ applied to~$v_\bpi$. Such a number is well-defined 
since~$V(\bpi)=\whbu_q^0(\ll)v_\bpi$ by~\propref{REP30}.
Set~$n(0)=-\infty$ with the convention that~$-\infty<n$ for all~$n\in\bz$.
Finally, set~$n(\bpi)=\max\{n(v)\,:\, v\in V(\bpi)\}$.
\begin{lem}
We have, for all~$v\in V(\bpi)$,
\begin{alignat*}{2}&n(E_iv)\le n(v),& \qquad &i\in I\\
&n(E_0v)\le n(v)+1,\\
&n(F_i v)\le n(v)+\dim V(\bpi)_{\lambda_\bpi-\alpha_i},&& i\in I\\
&n(F_0v)\le n(v)-1.
\end{alignat*}
\end{lem}
\begin{pf} 
The first two statements are obvious. 
For the next two, observe that since $V(\bpi)$ is spanned by vectors of the form~$Xv_\bpi$
where~$X$ is a monomial in the~$E_i$, $i\in\wh I$, it suffices by the relations in $\whbu_q$  to prove the assertion
for $v=v_\bpi$. If $i=0$, then $F_0v_\bpi=0$ and there is nothing to prove. So assume that $i\ne 0$, and that
$F_iv_\bpi\ne 0$ (if $F_iv_\bpi=0$, there is again nothing to prove).
Then,  by~\lemref{REP20}, $F_i v_{\bpi}$ is contained 
in the linear span of the~$x_{i,s}^-v_\bpi$, $s=1,\dots,\dim V(\bpi)_{\lambda_\bpi-
\alpha_i}$. 
\end{pf}

\subl{REP40}
Let~$m=m(\bpi)\in\bn^+$ be maximal such that~$\bpi\in (\bc(q)[u^m])^\ell$.
Then~$\bpi$ can be written uniquely as~$\bpi^0\bpi^0_\zeta\cdots
\bpi^0_{\zeta^{m-1}}$ where~$\bpi^0$ is an $\ell$-tuple of polynomials
with constant term~$1$, $\zeta$ is an~$m$th primitive root of
unity and the product is taken component-wise. 
By~\cite{C}, $V(\bpi)\cong V(\bpi^0)\tensor\cdots\tensor
V(\bpi^0_{\zeta^{m-1}})$. It follows
that~$\dim V(\bpi)_{\lambda_{\bpi}-\alpha_i}=m
\dim V(\bpi^0)_{\lambda_{\bpi^0}-\alpha_i}$ for all~$i\in I$.

Denote by~$\tau_{\bpi}$ the unique isomorphism of~$\whbu_q'$-modules
$$
V(\bpi^0)\tensor\cdots\tensor V(\bpi^0_{\zeta^{m-1}})
\longrightarrow
V(\bpi^0_{\zeta^{m-1}})\tensor V(\bpi^0)\tensor\cdots \tensor
V(\bpi^0_{\zeta^{m-2}})
$$
which sends~$v_{\bpi}=v_{\bpi^0}\tensor\cdots v_{\bpi^0_{\zeta^{m-1}}}$
to the corresponding permuted tensor product of highest weight vectors.
Set~$\eta_{\bpi}=(\phi_\zeta^*)^{\tensor m}\circ\tau_{\bpi}$,
where~$\phi_\zeta^*$ is the pull-back by the automorphism~$\phi_\zeta$
of~$\whbu_q$. Then~$\eta_{\bpi}(xv)=\zeta^{-\deg x} x\eta_{\bpi}(v)$
and~$\eta_\bpi(v_\bpi)=v_\bpi$, whence~$\eta_{\bpi}^m=\id$ and
$$
V(\bpi)=\bigoplus_{k=0}^{m-1} V(\bpi)^{(k)},\qquad\text{where~$V(\bpi)^{(k)}=
\{v\in V(\bpi)\,:\, \eta_\bpi(v)=\zeta^k v\}$}.
$$
Notice also that, since~$\deg K_i=0$, $\eta_{\bpi}$ preserves weight spaces
of~$V(\bpi)$. 
\begin{lem}
Let~$v\in V(\bpi)^{(k)}$ and suppose that~$v=\sum_{s=1}^N X_s
v_\bpi$ with~$X_s\in\whbu_q'$ homogeneous. Then $X_s v_\bpi\not=0$
only if~$\deg X_s+k=0\pmod m$. In particular, $n(v)+k=0\pmod m$.
\end{lem}
\begin{pf}
This is immediate since $V(\bpi)=\bigoplus_{r=0}^{m-1} V(\bpi)^{(r)}$ and $X_s v_\bpi\in
V(\bpi)^{(l)}$ where $l=-\deg X_s\pmod m$.
\end{pf}

\subl{REP50}
Let~$L(V(\bpi))=V(\bpi)\tensor_{\bc(q)} \bc(q)[t,t^{-1}]$ be the loop 
space of $V(\bpi)$.
Define the~$\whbu_q$-module structure on~$L(V(\bpi))$ by
$$
x(v\tensor t^n)=xv \tensor t^{n+\deg x},\quad D(v\tensor t^n)=q^n v\tensor 
t^n,\quad
C^{\pm1/2} (v\tensor t^n)=v\tensor t^n,
$$
for all~$x\in\whbu_q$ homogeneous, $v\in V(\bpi)$ and~$n\in\bz$.
Henceforth we write~$v t^n$ for the element~$v\tensor t^n$, $v\in V(\bpi)$, 
$n\in\bz$ of~$L(V(\bpi))$.

Let~$m=m(\bpi)$. By~\cite{CG}, $L(V(\bpi))$ is a direct sum of simple
submodules~$L^r (V(\bpi))$, $r=0,\dots,m-1$ where~$L^r(V(\bpi))=
\whbu_q(v_\bpi t^r)=\whbu_q(v_{\bpi}^*t^r)$.

Define~$\wh\eta_{\bpi}(vt^r)=\zeta^r \eta_{\bpi}(v)t^r$. 
Then by~\cite{CG}, $\wh\eta_{\bpi}\in\operatorname{End}_{\whbu_q} L(V(\bpi))$ and the 
simple submodule~$L^s(V(\bpi))$
is just the eigenspace of~$\wh\eta_{\bpi}$ corresponding to the 
eigenvalue~$\zeta^s$. 
\begin{lem}  For all $s=0,\dots ,m-1$, the $\bc(q)$-subspace 
$\wh\bu_q^+(v_\bpi t^s)$ is spanned by elements of the form $vt^{s+n(v)+k}$, 
$v\in V(\bpi)$, $k\in\bn$.\end{lem}
\begin{pf}
The statement follows immediately from~\corref{REP30} 
and from the defintion of $n(v)$ in~\ref{REP35}.
\end{pf}

\section{Irreducibility of~$X(\Lambda)\tensor V(\bpi)$}\label{SF}

In this section we prove the following
\begin{thm}\label{thmA}
Let~$\Lambda\in\wh P$ and let~$V(\bpi)$ be a finite dimensional
simple $\whbu_q'$-module corresponding to an $\ell$-tuple~$\bpi$ of
polynomials in one variable with constant term~$1$.
Then~$X(\Lambda)\tensor V(\bpi)$ is a simple $\whbu_q'$-module.
\end{thm}
\noindent
This result is a quantum version of~\cite[Theorem~4.2]{CPnew}

\subl{SF20}
By~\propref{P25}, $\whbu_q^r(\gg)$, $r\ge 0$ and $\whbu_q^r(\ll)$, $r > 0$
are contained in~$\whbu_q^+$ which is in turn contained
in the
Hopf subalgebra~$\wh\Cu_q^+$ of~$\whbu_q'$.
\begin{prop}
Let~$\Lambda\in\wh P$. 
\begin{enumerit}
\item As  $\whbu_q'$-module, we have
$$M(\Lambda)\otimes V(\bpi)=\whbu_q'(m_\Lambda\otimes v_\bpi)= \whbu_q'(m_\Lambda\otimes v_\bpi^*),$$
and 
$$X(\Lambda)\otimes V(\bpi)=\whbu_q'(v_\Lambda\otimes v_\bpi)= \whbu_q'(v_\Lambda\otimes v_\bpi^*).$$

\item As  $\whbu_q$-modules, we have
$$M(\Lambda)\tensor L^s(V(\bpi))=\sum_{n\in\bz} \whbu_q(m_\Lambda\otimes v_\bpi t^{mn+s})=\sum_{n\in\bz} \whbu_q(m_\Lambda\otimes v_\bpi^* t^{mn+s}),$$
and
 $$X(\Lambda)\tensor L^s(V(\bpi))=\sum_{n\in\bz} \whbu_q(v_\Lambda\otimes v_\bpi t^{mn+s})=\sum_{n\in\bz} \whbu_q(v_\Lambda\otimes v_\bpi^* t^{mn+s}),$$ for all $s=0,\cdots ,m-1$.
\end{enumerit}

\end{prop}
\begin{pf} 
The argument repeats that of the proof of~\cite[Lemma~2.1]{CPnew} and
is included here for the reader's convenience.

Observe first that, since~$\wh\Cu_q^+$ is a Hopf subalgebra of~$\whbu_q'$,
we have by~\lemref{REP30}
$$\wh\Cu_q^+(m_\Lambda\tensor v_{\bpi})
=m_\Lambda\tensor \whbu_q^+ v_{\bpi}=m_\Lambda\tensor V(\bpi).$$
 We prove by induction on~$\Ht\wh\gamma$ that
$$
M(\Lambda)_{\Lambda-\wh\gamma}\tensor V(\bpi)\subset
\whbu_q'(m_\Lambda\tensor V(\bpi)).
$$
If~$\Ht\wh\gamma=0$ then there is nothing to prove. Suppose that~$\Ht
\wh\gamma=1$ that is~$\wh\gamma=\alpha_i$ for some $i\in\wh I$.
Since~$M(\Lambda)_{\Lambda-\alpha_i}$ is spanned by~$F_i m_{\Lambda}$,
we have~$$F_i(m_\Lambda\tensor V(\bpi))=F_i m_{\Lambda}\tensor
V(\bpi)+m_{\Lambda}\tensor V(\bpi),$$ whence~$$F_i m_{\Lambda}
\tensor V(\bpi)\subset \whbu_q'(m_\Lambda\tensor V(\bpi)).$$ The inductive
step is proved similarly.

Thus, $M(\Lambda)\tensor V(\bpi)\subset \whbu_q'(m_\Lambda\tensor
v_{\bpi})$ we conclude that~$M(\Lambda)\tensor V(\bpi)$ is generated by~$m_\Lambda
\tensor v_{\bpi}$.
To see that it is also generated
by $m_\Lambda\otimes v_\bpi^*$ one proceeds as
above using an observation that, by~\propref{REP30},
$\whbu_q^0(\gg)v_{\bpi}^*=V(\bpi)$ and that $\whbu_q^0(\gg)m_\Lambda=0$. This proves~(i) for the 
modules $M(\Lambda)\otimes V(\bpi)$ and hence for the quotient module~$X(\Lambda)\otimes V(\bpi)$.

The proof of (ii) is similar. 
To see that induction starts, notice that by~\corref{REP30} and~\lemref{P30},
we have
$$
\sum_{n\in\bz} \whbu_q(m_\Lambda\otimes  v_\bpi t^{mn+s})=\sum_{n\in\bz}\whbu_q^-(m_\Lambda\otimes \whbu_q^+( v_\bpi t^{mn+s})=\sum_{n\in\bz}\whbu_q(m_\Lambda\otimes L^s(V(\bpi)).
$$
The inductive step is now completed as before.
\end{pf}

\subl{SF30}
Let~$\Lambda\in\wh P^+$. 
Recall from~\eqref{REP10.20}
 that when we regard~$X(\Lambda)$ as a~$\whbu_q'$-module, 
any weight vector $u\in X(\Lambda)$ of weight $\Lambda-\wh\gamma$ can be  
written uniquely as a sum ~$u=\sum_k u_{k}$ of linearly independent
elements  
$u_k\in X(\Lambda)_{\Lambda-\gamma_k-n_k\alpha_0}$,
where $\gamma_k\in Q^+$, $n_k\in\bn$ and $\gamma_k+n_k\alpha_0=\wh\gamma$. 
Denote by $\deg u$ the maximal value of the $n_k$. 

Given a weight vector~$v\in V(\bpi)$, set~$\Ht_\bpi(v):=
\Ht(\lambda_\bpi-\wt v)$.
Let $w\in X(\Lambda)\otimes V(\bpi)$ be a weight vector and write $$w=\sum_{k=1}^ru_k\tensor v_k$$ 
where the $u_k$ are linearly independent weight vectors in 
$X(\Lambda)$ and the $v_k$ are weight vectors in $V(\bpi)$. 
Fix  an integer $j_0(w)=j_0$,  $1\le j_0\le r$ such that the following two conditions hold
\begin{align}\lbl{1}
&\deg u_{j_0}\ge \deg u_j,\qquad \forall\, 1\le j\le r,
\\
\lbl{2}  
&\deg u_{j_0}=\deg u_j
\implies \Ht_\bpi(v_{j_0})\ge \Ht_\bpi(v_j).
\end{align}
\begin{prop} Let $w=\sum_{k=1}^r u_k\otimes v_k\in X(\Lambda)\otimes V(\bpi)$ be a weight
vector and let $j_0=j_0(w)$ 
be as above.
\begin{enumerit}
\item Assume that $v_{j_0}\notin \bc(q) v_\bpi$. Then there exists $i\in I$, $s\ge 0$
such that 
\begin{equation}\lbl{3}
0\ne x_{i,s}^+w=\sum_{\substack{j\,:\,\deg u_j=\deg u_{j_0},\\
\Ht_\bpi(v_j)=\Ht_\bpi(v_{j_0})}}q_i^{\wt(u_j\alpha_i^\vee)}(u_j\otimes x_{i,s}^+ v_j) + S,
\end{equation}
where $S$ is a sum of terms of the form $u_j'\otimes v_j'$ with either
$\deg u_j'<\deg u_{j_0}$ or $\deg u_j'=\deg u_{j_0}$
and $\Ht_\bpi(v_j')<\Ht_\bpi(x_{i,s}^+v_{j_0})$.

\item Suppose that $\deg u_{j_0}=0$. Then, for all  $s\gg 0$ there exist
$x\in\whbu_q^s(\gg)$ such that $xw=v_\Lambda\otimes v_\bpi$ and $y\in\whbu_q^s(\ll)$ such that 
$yw=v_\Lambda\otimes v_\bpi^*$.

\item Suppose that $\deg u_{j_0}=N$.  Then, there exists $s>0$ and 
$x\in\whbu_q^s(\gg)$  such that $xw=v_\Lambda\otimes v_\bpi$ and an element
$y\in\whbu_q^s(\ll)$ such that $yw=v_\Lambda\otimes v_\bpi^*$.
\end{enumerit}
\end{prop}
\begin{pf} 
By~\ref{REP10} and~\propref{REP30}
there exist $i\in I$ and $s\in\bn^+$ such that $x_{i,s}^+u_j=0$ for all $j$ and  
$x_{i,s}^+v_{j_0}\ne 0$. 
Observe that by~\eqref{P30.10b}, 
$$
x_{i,s}^+(u_j \tensor v_j)= q_i^{\wt u_j(\alpha_j^\vee)}
(u_j\otimes x_{i,s}^+v_j)+\sum_k u_k'\tensor v_k'
$$
where~$\Ht_\bpi(v_k')< \Ht_\bpi(v_j)$ and $\deg u_k'<\deg u_j$.
It follows that we can write $x_{i,s}^+w$ as in~\loceqref{3}. 
Notice that the  term $u_{j_0}\otimes v_{j_0}$ occurs with a non-zero coefficient on the right 
hand side of~\loceqref{3} and  is clearly linearly independent from the other terms in this equation. 
Hence $x_{i,s}^+w\ne 0$ and (i) is proved. 

To prove (ii), notice that if $\deg u_{j_0}=0$, then $\deg u_j=0$ for all $j$ and hence $xu_j=0$ 
for all $x\in\whbu_q^+$ which are homogenous of positive degree. 
It follows from~\eqref{P30.10b} that for all $x\in\whbu_q^1(\gg)$ 
we have $$xw=\sum_k u_k\otimes xv_k.$$ 
Choose $k$ such that $\Ht_\bpi(v_k)$ is maximal and $x\in \whbu_q^1(\gg)$ such 
that $xv_k=v_\bpi$. Then for all $j$ we have $xv_j=a_jv_\bpi$ for some $a_j\in\bc(q)$. 
It follows that $xw=u\otimes v_\bpi$ for some $u\in X(\Lambda)$ and $\deg u=0$. Since $\deg u=0$ 
and~$X(\Lambda)$ is irreducible 
it follows that $u\in U_q(\lie g)v_\Lambda$ and hence there exists $x'\in U_q(\lie g)\cap \whbu_q^+$ 
such that $x'u=v_\Lambda$ and so we get $x'xw=v_\Lambda\otimes v_\bpi$. 
Furthermore, there exists $y\in\whbu_q^1(\ll)$ such that $yv_\bpi=v_\bpi^*$. It follows 
that $y(v_\Lambda\otimes v_\bpi)=v_\Lambda\otimes v_{\bpi^*}$ which completes the proof of~(ii).

We prove (iii) by induction on $N$. Notice that~(ii) proves that induction starts..
Consider first the case when $v_{j_0}=a v_\bpi$ for some~$a\in\bc(q)^\times$.
Then we can write  $$w=u\otimes v_\bpi +w'$$ where 
$w'=\sum_{j\,:\,\deg u_j<N} u_j\otimes v_j.$  
Choose $s>0$ so that $x_{i,s}^+u=0$ for all $i\in I$. By~\eqref{P30.10b},
 $x^+_{i,s}(u\otimes v_\bpi)=0$. The induction hypothesis applies to $w'$ and 
so there exists $x\in\whbu_q^s(\gg)$ such that $xw'=v_\Lambda\otimes v_\bpi$. 
It follows that $xw=v_\Lambda\otimes v_\bpi$ and we are done.

Suppose then that $\Ht_\bpi(v_{j_0})=M$ and that~(iii) is established for $\Ht_\bpi(v_{j_0})<M$. 
By part (i) there exist $i\in I$ and $s>0$ such that $w'=x_{i,s}^+w\ne 0$. Furthermore 
write $w'=\sum_j u_j'\otimes v_j'$ and set $j_0'=j_0(w')$. Observe that
$\deg u'_{j_0'}=\deg u_{j_0}$ and $\Ht_\bpi (v'_{j_0'})=M-1$. Hence the induction hypothesis on $M$ 
applies and we conclude that there exists $x'\in\whbu_q^{s'}(\gg)$, $s'>s$ 
with $x(x_{i,s}^+) w=v_\Lambda\otimes v_\bpi$. 
\end{pf}
\begin{cor}
Let $W$ be a non-zero submodule of $X(\Lambda)\tensor V(\bpi)$ with~$\Lambda\in\wh P$
dominant. Then $W$ contains both $v_\Lambda\tensor v_\bpi$ and $v_\Lambda\tensor
v_\bpi^*$.
\end{cor}
\noindent
\thmref{thmA} follows immediately from the above Corollary and~\propref{SF20}.

\section{Filtrations of $M(\Lambda)\otimes L(V(\bpi))$ and $X(\Lambda)\otimes L(V(\bpi))$}\label{FIL}

\subl{FIL30}
Let~$M$ be a $\whbu_q$-module.
\begin{defn}
We call a collection~$\{\Cf_n\}_{n\in\bz}$ of
$\whbu_q$-submodules of~$M$ a decreasing 
$\bz$-filtration of~$M$ if $M=\sum_{n\in\bz} \Cf_n$ and
$\Cf_n\supseteq \Cf_{n+1}$ for all~$n\in\bz$. We say 
that the filtration $\{\Cf_n\}_{n\in\bz}$ is {\em strictly decreasing} if
$\Cf_n\ne \Cf_{n+1}$ for all $n\in\bz$ and is {\em trivial} if
$\Cf_m=\Cf_n$ for all~$m,n\in\bz$.
\end{defn}

In this section we prove that for $s=0,\cdots ,m-1$, the modules $M(\Lambda)\otimes L^s(V(\bpi))$
and $X(\Lambda)\otimes L^s(V(\bpi))$ 
admit $\bz$-filtrations $\Cm_n^{(s)}$ (respectively, $\Cx_n^{(s)}$), $n\in\bz$,
whose successive quotients are in the category $\Co$ and are 
isomorphic as $\whbu_q'$-modules. 
We  prove that $\Cm_n\supsetneq\Cm_{n+1}$ for all $n\in\bz$ and that
  $\bigcap_{n\in\bz}\Cm_n^{(s)}=0$
We also show that  the filtration $\Cx_n$, $n\in\bz$
 is either trivial or strictly decreasing. In the first case we prove that this implies that $X(\Lambda)\otimes L^s(V(\bpi))$ is irreducible and  
in the second  case we prove that  $\bigcap_{n\in\bz}\Cx_n^{(s)}=0$.

\subl{FIL35}
Set~$m=m(\bpi)$. 
\begin{prop}
Let~$\Lambda\in\wh P$. Given~$n\in\bz$,
let~$\Cm_n$ be the~$\whbu_q$-submodule of~$M(\Lambda)\tensor L(V(\bpi))$
generated by the vectors~$m_\Lambda\tensor v_{\bpi}t^{m n+s}$, $s=0,\dots,m-1$.
Then the modules~$\Cm_n$ form a $\bz$-filtration of~$M(\Lambda)
\tensor L(V(\bpi))$. Moreover, for all~$n\in\bz$, the modules $\Cm_n/\Cm_{n+1}$
are in the category~$\Co$ and are isomorphic as~$\whbu_q'$-modules. 

Further, if~$\Cx_n$ is the submodule of ~$X(\Lambda)\tensor L(V(\bpi))$
generated by the vectors~$v_\Lambda\tensor v_{\bpi}t^{m n+s}$, $s=0,\dots,m-1$,
then the modules~$\Cx_n$ form a $\bz$-filtration of~$X(\Lambda)
\tensor L(V(\bpi))$. Moreover, for all~$n\in\bz$ the modules $\Cx_n/\Cx_{n+1}$
are in the category~$\Co$ and are isomorphic as~$\whbu_q'$-modules. 
\end{prop}
\begin{pf}
We prove only the statement for the Verma modules, the proof of the one for~$X(\Lambda)$
being similar. Let~$\bpi=(\pi_i(u))_{i\in I}$, where~$\pi_i(u)=\sum_{k} 
\pi_{i,k} u^k\in\bc(q)[u]$. 
By the choice of~$m$, there exists~$i\in I$ such that~$\pi_{i,m}\not=0$. Then
$P_{i,m} v_{\bpi}=\pi_{i,m} v_\bpi$ and so
$P_{i,m}(m_\Lambda
\tensor v_{\bpi}t^{mn+s})=m_\Lambda\tensor (\pi_{i,m} v_{\bpi})t^{m(n+1)+s}$ by~\eqref{P30.10d}.
Therefore, $\Cm_n \supseteq \Cm_{n+1}$.
 
To show that $\Cm_n/\Cm_{n+1}$ is in the category $\cal{O}$, it
suffices to prove that the subspaces $\wh\Cu_q^+(m_\Lambda\tensor
v_\bpi t^{m n+s})$, $s=0,\dots,m-1$ of~$\Cm_n$ are finite-dimensional
modulo~$\Cm_{n+1}$.
Equivalently, it is sufficient to prove that the subspace
$\whbu_q^+(v_\bpi t^{m n})$ is
finite-dimensional modulo the subspace $\whbu_q^+(v_\bpi t^{m(n+1)})$. 
Now, by~\propref{P25},
$$
\whbu_q^+(v_\bpi t^{m n})\subset
\whbu_q^0(\ll)\whbu_q^0(0)(v_\bpi t^{m n})
=\whbu_q^0(\ll) v_{\bpi}t^{m n}\pmod{\whbu_q^+ v_\bpi t^{m(n+1)}}
$$
since~$P_{i,r}v_{\bpi}=0$ unless~$r$ is divisible by~$m$.
Since~$V(\bpi)$ is
finite dimensional, by~\lemref{REP30} 
there exist homogeneous~$X_1,\dots,X_N\in\whbu_q^r(\ll)$
for some~$r\ge 0$ such that~$X_1 v_{\bpi},\dots, X_N v_{\bpi}$
form a basis of~$V(\bpi)$.
Let~$x\in\whbu_q^0(\ll)$. Then there exist~$a_j\in\bc(q)$, $j=1,\dots,N$
such that
$$
x v_{\bpi}=\sum_{j=1}^N a_j X_j v_{\bpi}.
$$
We may assume that~$x$ is homogeneous of degree~$k$. 
Then~$x v_\bpi\in V(\bpi)^{(-k)}$ and so~$a_j=0$ unless~$\deg X_j=k\pmod m$ by~\lemref{REP40}. Then
we can write
$$
x(v_\bpi t^{m n})=(xv_\bpi)t^{m n+k}=\Big(\sum_{j=1}^N a_j X_j
v_\bpi\Big)t^{m n+k}= \sum_{j=1}^N a_j X_j(v_\bpi
t^{m n+k-\deg X_j}),
$$
the only non-zero terms being those with~$\deg X_j=k\pmod m$.
It follows that $x(v_\bpi t^{m n})=0\pmod{\whbu_q^+(v_\bpi t^{m(n+1)})}$
if~$k$ is sufficiently large.
 Therefore, the dimension
of $\whbu_q^+(v_\bpi t^{mn})\pmod{\whbu_q^+(v_\bpi t^{m(n+1)})}$ is bounded above 
by the dimension of the subspace of~$\whbu_q^0(\ll)$ spanned by
homogeneous elements whose degree does not exceed~$\max_j\{\deg X_j\}$.
Evidently, such a subspace of~$\whbu_q^0(\ll)$ is finite-dimensional.

To prove that $\Cm_n/\Cm_{n+1}\cong \Cm_{n-1}/\Cm_n$ as a~$\whbu_q'$-module
for all $n$, consider the map
\begin{alignat*}{2}
M(\Lambda)\otimes L(V(\bpi))&\longrightarrow && M(\Lambda)\otimes
L(V(\bpi))\\
v\tensor w t^k &\mapsto && v\tensor w t^{k+m},
\end{alignat*}
for all~$v\in M(\Lambda)$, $w\in V(\bpi)$ and~$k\in\bz$.
This is obviously a map of $\whbu_q'$-modules (but, of course, not 
$\whbu_q$-modules) which 
 takes $\Cm_n$
isomorphically onto $\Cm_{n+1}$. Moreover, this operation corresponds
to tensoring~$\Cm_n$
 with the 1-dimensional highest weight integrable module~$X(m\delta)$.
Thus  we have $\Cm_n\cong \Cm_{n+1}\tensor X(-m\delta)$ and so in fact
$\Cm_n/\Cm_{n+1}$ is isomorphic to~$(\Cm_{n-1}\tensor X(m\delta))/(\Cm_n
\tensor X(m\delta))$ as a $\whbu_q'$-module.
\end{pf}

\subl{FIL31}
For $s=0,\cdots ,m-1$, let~$\Cm_n^{(s)}$
 (resp. $\Cx_n^{(s)}$ be the $\whbu-q$-submodule
of  $M(\Lambda)\otimes L(V(\bpi))$ (respectively, of $X(\Lambda)\otimes L(V(\bpi))$)  
generated by~$m_\Lambda\tensor v_\bpi t^{m n+s}$ (respectively, by  $v_\Lambda\otimes v_\bpi t^{mn+s}$). 
\begin{lem} For all $n\in\bz$, we have 
$$\Cm_n=\bigoplus_{s=0}^{m-1}\Cm_n^{(s)},$$
and
$$\Cm_n^{(s)}=\Cm_n\cap (M(\Lambda)\otimes L^s(V(\bpi))).$$
Further, the $\Cm_n^{(s)}$, $n\in\bz$ form a decreasing filtration of 
$M(\Lambda)\otimes L^s(V(\bpi))$, $s=0,\cdots ,m-1$. Analogous statements hold for $\Cx_n^{(s)}$.
\end{lem} 
\begin{pf}
This follows immediately from the trivial observation that for any $\whbu_q$-module $M$ we have
$M\otimes L(V(\bpi))=\bigoplus_{s=0}^{m-1} M\otimes L^s(V(\bpi))$.
\end{pf}

\subl{FIL40}

Fix an ordered basis of $V(\bpi)$ of
weight vectors $v_0=v_\bpi$, $v_1,\dots ,v_N=v_\bpi^*$ such that
$\Ht_\bpi(v_i) \le \Ht_\bpi(v_{i+1})$
for all~$i=0,\dots,N-1$.
Furthermore, we may assume, without loss of generality,
that~$v_j\in V(\bpi)^{(k_j)}$ for some~$0\le k_j\le m-1$. 
It is clear  that~$F_i v_j$ is a linear combination of~$v_{j'}$ with~$j'>j$
if~$i\in I$ and with~$j'<j$ if~$i=0$. Let $\tilde\Cm_n^{(s)}$ be the 
 $\wh\Cu_q^-$-submodule of $M(\Lambda)\otimes L(V(\bpi))$ 
generated by the set~$\{m_\Lambda\otimes v_j t^{r}\,:
\, r \ge m n+s,\, j=0,\dots, N,\,r=s-k_j\pmod m\}$ and set $\sum_{s=0}^{m-1}
\tilde\Cm_n^{(s)}=\Cm_n$. Similarly,
let $\tilde\Cx_n^{(s)}$ be the  $\wh\Cu_q^-$-submodule of~$X(\Lambda)\tensor L^s(V(\bpi))$ generated
by the set~$\{v_\Lambda\tensor v_j t^r\,:\, r\ge m n+s,\, j=0,\dots,N,\, r=s-k_j\pmod m\}$. 
\begin{lem} 
For~$0\le s\le m-1$, and $n\in\bz$, we have
$$\Cm_n^{(s)}\subset\tilde\Cm_n^{(s)},\qquad\Cx_n^{(s)}\subset\tilde\Cx_n^{(s)}.$$
\end{lem} 
\begin{pf}
Immediate.
\end{pf}

\subl{FIL50} 
The following proposition plays a crucial role in the reminder of the paper.
\begin{prop}
\begin{enumerit}
\item Let $v\in V(\bpi)$, $n\in\bn$ and
suppose that there exist elements $X_{j,r}\in(\wh\bu_q^-)_+$,
$r\ge 0$, $j=0,\dots ,N$ such that in $M(\Lambda)\otimes L(V(\bpi))$ we have,
\begin{equation}\lbl{11}m_\Lambda\otimes vt^n=\sum_{j=0}^N\sum_{r\ge n}X_{j,r}
(m_\Lambda\otimes v_jt^r).
\end{equation}
Then $v=0$.
\item Let $w\in M(\Lambda)\otimes L(V(\bpi))$, $w\not=0$. Then there 
exists $n\in\bz$ such that  $w\notin \tilde\Cm_{n+1}$.
\item  Let $v\in V(\bpi)$ and
suppose that there exist elements $X_{j,r}\in(\wh\bu_q^-)_+$,
$r\ge 0$, $j=0,\dots ,N$ such that in $X(\Lambda)\otimes L(V(\bpi))$ we have,
\begin{equation}\lbl{12}v_\Lambda\otimes vt^n=\sum_{j=0}^N\sum_{r\ge n}X_{j,r}
(m_\Lambda\otimes v_jt^r).
\end{equation}
Let $R$ be the maximal~$r\ge n$ such that there exists $0\le j\le N$ with $X_{j,r}
( v_\Lambda\otimes v_{j}t^{r})\ne 0$ and let $j_0$ be the minimal~$j$
such that~$X_{j,R}(v_\Lambda\tensor v_j t^{R})\not=0$.
Then $X_{j_0,R}\in\Ann_{\wh\bu_q^-} v_\Lambda$. 
\end{enumerit} 
\end{prop}
\begin{pf} 
To prove~(i), suppose for a contradiction that $v\ne 0$.
Let $R$ be the maximal~$r\ge n$ such that there exists $j$  with $X_{j,r}(m_\Lambda\otimes v_j t^r)\ne 0$ and 
let $j_0$ be the minimal~$0\le j\le N$ such that~$X_{j,R}(m_\Lambda\tensor v_{j}t^R)\not=0$.
Then $X_{j_0,R}
(m_\Lambda\tensor v_{j_0} t^{R})$ contains a term
$c(X_{j_0,R}m_\Lambda)\tensor v_{j_0} t^{R}$ 
for some~$c\in\bc(q)^\times$. Since
$\Delta(\wh\bu_q^-)\subset\wh\Cu_q^-\otimes \wh\Cu_q^-$, it follows that 
all other elements in~\loceqref{11} 
are terms of the form $ m'\otimes v_{j'}t^r$ where either $r<R$ or $r=R$
and $j'>j_0$. If $R>n$, then 
this forces~$X_{j_0,R}\in \Ann_{\whbu_q^-}m_\Lambda$
and hence $X_{j_0,R}=0$ which is a contradiction.
If $R=n$, then~\loceqref{11} reduces to
$$
m_\Lambda\otimes vt^n=
\sum_{j=0}^N X_{j,n}(m_\Lambda\otimes v_j t^n).
$$ 
Let  $k=\#\{j: X_{j,n}(m_\Lambda\otimes v_jt^n)\ne 0\}$. If $k=0$ then we are done. 
Suppose that~$k=1$. Then $m_\Lambda\otimes vt^n=X_{j,n}(m_\Lambda\otimes v_jt^n)$ for some $0\le j\le N$.  
If  $X_{j,n}m_\Lambda\ne 0$ then 
$X_{j,n}(m_\Lambda\otimes vt^n)\in \bigoplus_{\gamma\in\wh{Q}^+\setminus\{0\}}M(\Lambda)_{\Lambda-\gamma}\otimes L(V(\bpi))$ 
which is clearly impossible. 
Hence $X_{j,n}m_\Lambda =0$ and we get a contradiction.
Suppose then that we have proved that either $k=0$ or $k\ge s$ for some~$s\in\bn^+$. 
If $k=s$ and $X_{j_r,n}(m_\Lambda\otimes vt^{j_r})\ne 0$ for some~$0\le j_1<\cdots< j_s\le N$.
If $X_{j_r,n}m_\Lambda\ne 0$ for any $j_r$, then again the right hand side of~\loceqref{11} is 
contained in~$\bigoplus_{\gamma\in\wh{Q}^+\setminus\{0\}}M(\Lambda)_{\Lambda-\gamma}\otimes L(V(\bpi))$
which is a contradiction. Thus $k=0$ or $k>s$. Since $V(\bpi)$ is finite-dimensional
it follows that $k=0$ and we are done.

To prove (ii), write $w=\sum_{j=1}^s m_j\otimes w_jt^{r_j}$, where $m_j\in M(\Lambda)$, $w_j\in V(\bpi)$ and $r_j\in\bz$. Let $n_0=\max\{r_j:1\le j\le s\}$ and suppose that $w\in\tilde\Cm_n$ for some $n>n_0$.
This means that we can write
$$w=\sum_{j=0}^N\sum_{r\ge n}X_{j,r}(m_\Lambda\otimes v_jt^r)$$
for some choice of $X_{j,r}\in\wh\bu_q^-$. But now arguing exactly as in the $R>n$ case of (i), we see that $w=0$ which is a contradiction.

The proof of (iii) is an obvious modification of the  argument in (i).
\end{pf}

\subl{FIL60}
\begin{prop}
The $\bz$-filtration $\Cm_n^{(s)}$ of~$M(\Lambda)\tensor L^s(V(\bpi))$, $s=0,\dots,m-1$
is strictly decreasing and $\bigcap_{n\in\bz}\Cm_n^{(s)}=0$.
\end{prop}
\begin{pf} 
In view of \propref{FIL35} for the first statement it is sufficient
 to prove that $\Cm_0^{(s)}\ne \Cm_1^{(s)}$.
 Assume for a contradiction that $\Cm_0^{(s)}=\Cm_1^{(s)}$. Then
 $m_\Lambda\otimes v_\bpi t^s\in\Cm_1^{(s)}$
 and hence it follows from~\lemref{FIL40} that
there exist~$X_{j,r}\in\whbu_q^-$ such that
\begin{equation}\lbl{10}
m_\Lambda\tensor v_0 t^s=\sum_{j=0}^N
\sum_{r \ge m+s} X_{j,r}(m_\Lambda\tensor v_j t^r). 
\end{equation}
Since~$X_{j,r}\in(\wh\bu_q^-)_+$  we get a contradiction by~\propref{FIL50}(i).

Furthermore,
let $w\in\bigcap_{n\in\bn}\Cm_n^{(s)}$. If $w\ne 0$, then by~\propref{FIL50}(ii)
we can choose $n_0\in\bz$ such that $w\notin \tilde\Cm_{n_0}^{(s)}$ contradicting $\Cm_{n_0}^{(s)}
\subset\tilde\Cm_{n_0}^{(s)}$.
\end{pf}

\subl{FIL65} To analyse the filtration on $X(\Lambda)\otimes L(V(\bpi))$, $\Lambda\in\wh P^+$, we need the 
following analogue of~\propref{SF30}. 
\begin{prop} 
Let $w=\sum_{k=1}^r u_k\otimes v_kt^{r_k}\in X(\Lambda)\otimes L(V(\bpi))$ be a weight vector
and let  $j_0=j_0(w)$ be the integer associated with the element 
$\sum_{k=1}^r u_k\otimes v_k\in X(\Lambda)\otimes V(\bpi)$ as in~\ref{SF30}.
\begin{enumerit}
\item Assume that $v_{j_0}\ne v_\bpi$.  There exists $i\in I$, $s\ge 0$ such that 
\begin{equation}\lbl{3}
0\ne x_{i,s}^+w=\sum_{\substack{j\,:\,\deg u_j=\deg u_{j_0},\\
\Ht_\bpi(v_j)=\Ht_\bpi(v_{j_0})}}q_i^{\wt(u_j\alpha_i^\vee)}
(u_j\otimes x_{i,s}^+ v_jt^{r_j+s}) + S,
\end{equation}
where $S$ is a sum of terms of the form $u_j'\otimes v_j't^{r_j'}$ where either
$\deg u_j'<\deg u_{j_0}$ or $\deg u_j'=\deg u_{j_0}$
and $\Ht_\bpi(v_j')<\Ht_\bpi(x_{i,s}^+v_{j_0})$.

\item Suppose that $\deg u_{j_0}=0$. Then $r_j=R$ for all $j$ for some $R\in\bn$.  
Furthermore, for all  $s\gg 0$ there exists 
$x\in\bu_q^s(\gg)$ and an integer $L$
 such that $xw=v_\Lambda\otimes v_\bpi t^L$ and an element 
$y\in\bu_q^s(\ll)$ and an integer $L'$ such that $yw=v_\Lambda\otimes v_\bpi^*t^{L'}$.

\item Suppose that $\deg u_{j_0}=N$.  Then there exists $s>0$ and 
$x\in\bu_q^s(\gg)$  and an integer $S$ such that $xw=v_\Lambda\otimes v_\bpi t^S$ and an element
$y\in\bu_q^s(\ll)$  and an integer $S'$ such that $yw=v_\Lambda\otimes v_\bpi^*t^{S'}$.
\end{enumerit}
\end{prop}
\begin{pf} 
The first statement in part (ii) is an obvious consequence of the fact that $w$ is a weight vector.  
The proposition is now  proved in exactly the same way as~\propref{SF30} and we omit the details.
\end{pf}
\begin{cor}
Let~$W$ be a non-zero submodule of~$X(\Lambda)\tensor L(V(\bpi))$ with~$\Lambda\in\wh P$
dominant. Then~$W$ contains
$v_\Lambda\tensor v_\bpi t^s$ for some~$s\in\bz$ and~$v_\Lambda\tensor v_\bpi^* t^{r}$
for some~$r\in\bz$.
\end{cor}

\subl{FIL66} We note the following consequence of~\corref{FIL65}
\begin{prop}
The $\bz$-filtration $\Cx_n^{(s)}$ of~$X(\Lambda)\tensor L^s(V(\bpi))$, 
$s=0,\dots,m-1$ is 
either
\begin{enumerit} 
\item  trivial
and $X(\Lambda)\otimes L^s(V(\bpi))$ is irreducible, or
\item
strictly decreasing and
$$\bigcap_{n\in\bz}\Cx_n^{(s)}=0.$$
\end{enumerit}
\end{prop}
\begin{pf}
Suppose that $\Cx_{n}^{(s)}=\Cx_{m}^{(s)}$ for some $m>n\in\bz$.
Then~$\Cx_{n}^{(s)}=\Cx_{n+1}^{(s)}$ and 
it follows from~\propref{FIL35} that $\Cx_m^{(s)}=\Cx_n^{(s)}$ for all
$m,n\in\bz$. This proves that the filtration is either trivial or strictly 
decreasing.

Let $W$ be a non-zero submodule of $X(\Lambda)\otimes L^s(V(\bpi))$. By~\corref{FIL65},
$v_\Lambda\otimes v_\bpi t^{m r+s}\in W$ for some~$r\in\bz$
and so $\Cx_{r}^{(s)}\subset W$. If the filtration is trivial, 
then this implies that $\Cx_n^{(s)}\subset W$ for all $n\in\bz$. It follows from~\propref{SF20}(ii) that 
$W=X(\Lambda)\otimes L^s(V(\bpi)$ and (i) is proved.

Suppose that the filtration is strictly decreasing and set $W=\bigcap_{n\in\bz}\Cx_n^{(s)}$. Suppose
that $W\not=0$. Then it follows from~\corref{FIL65} that $v_\Lambda\otimes v_\bpi t^{m r+s}\in W$ for some~$r\in\bz$ and 
hence $\Cx_{r}^{(s)}\subset W$. Then
$\Cx_r^{(s)}\subset\Cx^{(s)}_{r+1}$ and so~$\Cx_r^{(s)}=\Cx_{r+1}^{(s)}$ which is a contradiction whence~(ii).
\end{pf}

\subl{FIL70} The results of this section allow us to 
complete the tensor products $X(\Lambda)\otimes L(V(\bpi))$ and $M(\Lambda)\otimes L(V(\bpi))$. We restrict ourselves 
to the first case, the second one being similar.
Let~$\Lambda\in\wh P$ and suppose that the filtration~$\Cx_r^{(s)}$, $r\in\bz$ is strictly decreasing.
Let~$X(\Lambda)\wh\tensor L^s(V(\bpi))$ be the completion of~$X(\Lambda)\tensor L^s(V(\bpi))$ with
respect to the topology induced by the filtration~$\Cx_r^{(s)}$. It is well-known that
then there exists a canonical
map~$\phi_\Cx: X(\Lambda)\tensor L^s(V(\bpi))\longrightarrow X(\Lambda)\wh\tensor L^s(V(\bpi))$ and~$\ker\phi_\Cx=
\bigcap_{r\in\bz} \Cx_r^{(s)}=0$ by~\propref{FIL66}(ii). Therefore, $X(\Lambda)\tensor L^s(V(\bpi))$
embeds into the completion. On the other hand, 
$$
X(\Lambda)\wh\tensor L^s(V(\bpi))\cong \lim_{\longleftarrow} (X(\Lambda)\tensor L^s(V(\bpi)))/\Cx_r^{(s)}.
$$
Furthermore, let~$\wh\Cx_n^{(s)}$ be the completion of~$\Cx_n^{(s)}$,
$$
\wh\Cx_n^{(s)}=\lim_{\longleftarrow} \Cx_m^{(s)}/\Cx_n^{(s)}.
$$
Then~$\wh\Cx_n^{(s)}$ is a $\bz$-filtration on $X(\Lambda)\wh\tensor L(V(\bpi))$ and
$\Cx_n^{(s)}/\Cx_{n+1}^{(s)}\cong \wh\Cx_n^{(s)}/\wh\Cx_{n+1}^{(s)}$ 
and so the associated graded space of $X(\Lambda)\tensor L^s(V(\bpi))$ with respect to
the filtration~$\Cx_n^{(s)}$ is isomorphic to the associated graded space of~$X(\Lambda)\wh\tensor L(V(\bpi))$
with respect to the filtration~$\wh\Cx_n^{(s)}$.
One also has~$(X(\Lambda)\wh\tensor L(V(\bpi)))/\wh \Cx_n\cong (X(\Lambda)\tensor L(V(\bpi)))/\Cx_n$.

\section{An irreducibility criterion for~$X(\Lambda)\tensor L(V(\bpi))$}\label{IRR}

In this section we establish a sufficient condition for the simplicity of
the~$\whbu_q$-modules $X(\Lambda)\tensor L^s(V(\bpi))$,
$s=0,\dots,m(\bpi)-1$, $\Lambda\in \wh P^+$.
\begin{thm}\label{thmB}
Let $\Lambda\in\wh P^+$ and let $\bpi=(\pi_i(u))_{i\in I}$ be
an $\ell$-tuple of polynomials with constant term~$1$.
Suppose  that either
\begin{align*}
&(k(\bpi)+m(\bpi))(\Lambda,\delta) < (\Lambda+\lambda_\bpi,\alpha_i),
\\
\intertext{for some~$i\in I$ satisfies~$k(\bpi)=\dim V(\bpi)_{\lambda_{\bpi}-\alpha_i}$ or}
&k(\bpi)(\Lambda,\delta) < -(\Lambda+w_\circ\lambda_\bpi,\alpha_i),
\end{align*}
for some~$i\in I$ satisfying~$k(\bpi)=\dim V(\bpi)_{w_\circ\lambda_\bpi+
\alpha_i}$. For all $s=0,\cdots ,m-1$, the filtration $\Cx_n^s$, $n\in\bz$ of $X(\Lambda)\otimes L^s(v\bpi)$ is trivial  and hence $X(\Lambda)\otimes L^(V(\bpi))$ is an irreducible $\whbu_q$--module.

\end{thm}

\begin{pf} It suffices by~\propref{FIL66} to show that the filtration 
is trivial. Let $\bpi^0$ be as defined in~\ref{REP40}, and
set~$k=k(\bpi^0)$,  clearly  $k(\bpi)=m k$. Fix $i\in I$ so that $k(\bpi)=\dim V(\bpi)_{\Lambda_\bpi-\alpha_i}$.
By~\lemref{REP20},
$$
x_{i,m(k+1)}^- v_{\bpi}= \sum_{r=1}^{m k} a_r x_{i,r}^- v_{\bpi}
$$
for some~$a_r\in\bc(q)$, $r=1,\dots,m k$. Applying~\eqref{P30.10c} we see that
\begin{align*}
x_{i,m(k+1)}^- (v_\Lambda\tensor v_{\bpi}t^{mn+s})&=
v_\Lambda\tensor (x_{i,mk+1}^- v_{\bpi})t^{m(n+k+1)+s}\\
&=v_\Lambda\tensor \Big(\sum_{r=1}^{m k} a_r x_{i,r}^-
v_{\bpi}\Big)t^{m(n+k+1)+s}
\\
&=\sum_{r=1}^{m k} a_r x_{i,r}^- (v_\Lambda\tensor v_{\bpi}t^{m(n+k+1)+s-r}),
\end{align*}
which proves that $x_{i,m(k+1)}^- (v_\Lambda\tensor v_{\bpi}t^{mn+s})$
is contained in~$\Cx_{n+1}^{(s)}$. Consider 
the subalgebra of~$\whbu_q$ generated
by~$E=x_{i,m(k+1)}^-$, $F=x_{i,-m(k+1)}^+$ and~$K=C^{m(k+1)}K_i^{-1}$ which  is isomorphic to~$U_{q_i}(\lie{sl}_2)$ with standard
generators~$E,F,K^{\pm1}$.
 Note that in $\Cx_n^{(s)}/\Cx_{n+1}^{(s)}$ we have
 $$x_{i,m(k+1)}^- (v_\Lambda\tensor v_{\bpi}t^{mn+s})=0,\quad K(v_\Lambda\tensor v_{\bpi}t^{mn+s})=q_i^r(v_\Lambda\tensor v_{\bpi}t^{mn+s})=0,$$
where $d_ir= m(k+1)(\lambda,\delta)-(\lambda+
\lambda_{\bpi},\alpha_i) <0$.  
Since $\Cx_n^{(s)}/\Cx_{n+1}^{(s)}$ is an integrable module for $\whbu_q$ and hence.for this copy 
of $\bu_{q_i}(\lie{sl}_2)$.  But this forces
$$ v_\Lambda\tensor v_{\bpi}t^{mn+s}\in\Cx_{n+1}\cap X(\Lambda)\otimes L^s(V(\bpi)) =\Cx_{n+1}^{(s)},$$ 
and so  $\Cx_n^{(s)}=\Cx_{n+1}^{(s)}$.

The second assertion is proved similarly. 
We work with the elements $v_\Lambda\otimes v_\bpi^*t^{mn+s}$ and $x_{i,mk}^+$ and we omit the
details.
\end{pf}

\section{A reducibility criterion for~$X(\Lambda)\tensor L(V(\bpi))$}\label{RDC}

In this section we analyse the structure of~$X(\Lambda)\tensor
L(V(\bpi))$ and give a sufficient condition for 
the tensor product to be reducible.

\begin{thm}\label{thmC}
Let~$\Lambda\in\wh P$ be dominant and suppose that~$(\Lambda,\delta)\ge
(\Lambda+\lambda_{\bpi})(\theta^\vee)+m(\bpi)$
or, equivalently, $\Lambda(\alpha_0^\vee)\ge
\lambda_{\bpi}(\theta^\vee)+m(\bpi)$.
Then the modules~$X(\Lambda)\tensor L^s(V(\bpi))$, $s=0,\dots,m(\bpi)-1$
are reducible.
\end{thm}

The theorem is proved in the rst of the section.
\subl{RDC10}

\begin{lem}
Let~$\Lambda$ be a dominant weight.
Take~$v\in V(\bpi)$.
Then
$$
F_i^{\Lambda_i+1} (v_{\Lambda}\tensor
v t^r)=
\sum_{k=1}^{\Lambda_i+1}
c_k F_i^{\Lambda_i+1-k}(v_\Lambda\tensor (F_i^k v)t^{r-k\delta_{i,0}})
$$
where~$c_k\in\bc(q)$.
\end{lem}
\begin{pf}
Recall that~$F_i^{\Lambda_i+1}v_\Lambda=0$.
One has
$$
F_i^{\Lambda_i+1}(v_\Lambda\tensor v t^r)=
F_i^{\Lambda_i}(c_1 F_i v_\Lambda\tensor v t^r+ v_\Lambda\tensor
(F_i v)t^{r-\delta_{i,0}}),
$$
The second term has the desired form.
If~$\Lambda_i=0$ then the first term equals zero and we are  done. Otherwise,
we can write
$$
F_i^{\Lambda_i}(F_i v_\Lambda\tensor vt^r)=
F_i^{\Lambda_i-1}(c' F_i^2 v_\Lambda\tensor v t^r+
F_i v_\Lambda\tensor F_i vt^{r-\delta_{i,0}}).
$$
Clearly, the element  $F_i v_\Lambda\tensor F_i v t^{r-\delta_{i,0}}$
is a linear combination of~$F_i(v_\Lambda\tensor (F_i v) t^{r-\delta_{i,0}})$
and~$v_\Lambda\tensor (F_i^2 v) t^{r-2\delta_{i,0}}$ which are
both of the required form.

Now suppose that, for all~$k=1,\dots,s-1$, we can
write~$F_i^{\Lambda_i+1-k}(F_i^k v_\Lambda\tensor v t^r)$
as a linear combination of terms which have the required form
and~$F_i^{\Lambda_i-k}(F_i^{k+1}v_\Lambda\tensor vt^r)$. Then
$$
F_i^{\Lambda_i+1-s}( F_i^s v_\Lambda\tensor vt^r)=
F_i^{\Lambda_i-s}(F_i^{s+1} v_\Lambda\tensor v t^r+
F_i^{s} v_\Lambda\tensor F_i v t^{r-\delta_{i,0}}).
$$
Now,  the second term can be written as a linear combination of terms
which have the required form by the induction hypothesis. Hence we can repeat
 the process until we get to~$s=\Lambda_i$ in which
case~$F_i^{s+1}$ annihilates~$v_\Lambda$ and we are done.
\end{pf}

\subl{RDC15}
Set~$m=m(\bpi)$. 
By~\propref{FIL35}, 
the module~$X(\Lambda)\tensor L^s(V(\bpi))$ admits
a $\bz$-filtration $\Cx_n^{(s)}=\whbu_q (v_\Lambda\tensor v_{\bpi}t^{m n+s})$.
Let $v_0,\dots ,v_N$ be the basis of $V(\bpi)$ introduced in~\ref{FIL40}.
\begin{lem}  
Let  $v\in V(\bpi)$, $s,n\in\bn$.
Suppose that there exists $R\in\bn$ and 
$X_{j,r}\in \wh\bu_q^-$, $j=0,\dots ,N$, $n+s\le r\le R$
such that 
\begin{equation}\lbl{11}
v_\Lambda \tensor v t^n =\sum_{j=0}^N
\sum_{r= n+s}^{R} X_{j,r}(v_\Lambda\tensor v_j t^r).
\end{equation} 
Let $R_0$ be the minimal value of $R$ for which such  an
expression exists. Then
$R_0\le \lambda_\bpi(\theta^\vee)+n+s$.
\end{lem}

\begin{pf} 
Assume for a contradiction that $R_0> \lambda_\bpi(\theta^\vee)+n+s$. 
Let~$j_0$ be minimal such 
that~$X_{j_0,R_0}(v_\Lambda\tensor v_{j_0}t^{R_0})\not=0$.
By~\propref{FIL50}(ii)  
and~\propref{REP10}(ii) $X_{j_0,R_0}=\sum_{i\in\wh I}
y_i F_i^{\Lambda_i+1}$ for some~$y_i\in\whbu_q^-$.
If $i\in I$, then by~\lemref{RDC10}
$y_iF_i^{\Lambda_i+1}(v_\Lambda\otimes v_{j_0}t^{R_0})$ is a linear 
combination of the elements
$y_iF_i^{\Lambda_i-k+1}(v_\Lambda\tensor F_i^kv_{j_0} t^{r})$, 
$k=1,\cdots,\Lambda_i+1$.
But these terms are all of the form
$X_{j',R_0}'(v_\Lambda\otimes v_{j'}t^{R_0})$ with~$j'>j_0$ and~$X_{j',R_0}'\in\whbu_q^-$.
By~\lemref{RDC10} and!\lemref{REP28} we conclude that
$y_0F_0^{\Lambda_0+1}(v_\Lambda\tensor v_{j_0}t^{R_0})$ is   
a linear combination 
of terms of the form 
$y_0 F_0^{\Lambda_0+1-k}(v_\lambda\otimes F_0^k v_{j_0}t^{R_0-k})$ where
$1\le k\le \min\{\Lambda_0+1,\lambda_\bpi(\theta^\vee)\}$. 
Observe that $R_0-k\ge R_0-\lambda_\bpi(\theta^\vee)\ge n+s$.
Thus, we have obtained another expression of the form~\loceqref{11}
where either $R<R_0$ or $R=R_0$ and the minimal value of~$j$ such that~$X_{j,R_0}(v_\Lambda
\tensor v_j t^{R_0})\not=0$ is strictly greater than~$j_0$.
The former situation cannot occur by the choice of $R_0$. On the other hand,
since~$V(\bpi)$ is finite-dimensional, using the above argument repeatedly
we must eventually reach a stage where $n+s\le R<R_0$ which is a contradiction.
\end{pf}

\subl{RDC20}  \thmref{thmC} is an immediate consequence of the following
\begin{prop}
Suppose that~$\Lambda(\alpha_0^\vee)>\Lambda_{\bpi}(\theta^\vee)+m-1$.
Then the $\bz$-filtration 
$\Cx_n^{(s)}$ on~$X(\Lambda)\tensor L^s(V(\bpi))$ is strictly
decreasing.
\end{prop}
\begin{pf} 
Assume for a contradiction that $\Cx_0^{(s)}=\Cx_1^{(s)}$. 
Then, as in~\propref{FIL60}, we can write
\begin{equation}\lbl{10}
v_\Lambda\tensor v_0 t^s = \sum_{r=m+s}^{R}
\sum_{j=0}^N
X_{j,r}(v_\Lambda\tensor v_j t^r),
\end{equation}
for some $R\ge m+s$ and $X_{j,r}\in(\whbu_q^-)_+$ 
with $X_{j,R}(v_\Lambda\tensor v_j t^R)\ne 0$ for some $0\le j\le N$.

Assume that~$R$ is minimal such that the expression of the form~\loceqref{10}
exists. Then~$m+s\le R\le \lambda_\bpi(\theta^\vee)+m+s$ by~\lemref{RDC15}.
Furthermore, let~$j_0$ be the minimal value of $j$, such that
 $X_{j_0,R}(v_\Lambda\tensor v_j t^R)\ne 0$. Then by~\propref{FIL50}(ii)  
and~\propref{REP10}(ii), 
$X_{j,_0,R}=\sum_{i\in I} y_i F_i^{\Lambda_i+1}$ for some~$y_i\in\wh\bu_q^-$.
Furthermore, $X_{j_0,R}$ is of weight
$-(R-s)\delta+\wt v_0-\wt v_j\in-(R-s)\alpha_0+
Q^+$.
On the other hand, the weight of~$y_0 F_0^{\Lambda_0+1}$ is contained in the
set~$-(\Lambda_0+1)\alpha_0-\wh Q^+$. Since 
$R-s\le \lambda_\bpi(\theta^\vee)+m<\Lambda_0+1$ 
we conclude that~$y_0 F_0^{\Lambda_0+1}(v_\Lambda\tensor
v_j t^{R})=0$. It follows that~$X_{j_0,R}(v_\Lambda \tensor
v_j t^{R})=\sum_{i\in I}y_iF_i^{\Lambda_i+1}(v_\Lambda \tensor
v_j t^{R})$. Thus,  by~\lemref{RDC10} $X_{j_0,R}(v_\Lambda\tensor
v_j t^R)$ is a linear combination of terms of the form~$X_{j',R}(v_\Lambda\tensor v_{j'}t^{R})$
with~$j'>j$. Since~$V(\bpi)$ is finite-dimensional, repeating this process  we obtain
an expression of the form~\loceqref{10} with~$X_{j,R}(v_\Lambda\tensor v_j t^R)=0$
for all~$0\le j\le N$ which contradicts the minimality of~$R$.
\end{pf}

\section{ Structure of $\Cx_n/\Cx_{n+1}$ in some special cases.}\label{RN}

In this section we analyse the quotient modules $\Cx_n/\Cx_{n+1}$ in 
the special case when $\lie g$ is of 
type $A_\ell,B_\ell,C_\ell,D_\ell$ and
$V(\bpi)$ is isomorphic as a $U_q(\lie g)$-module to the natural 
representation of $U_q(\lie g)$.
\subl{RN10} Assume that the nodes of the Dynkin diagram of $\lie g$ are numbered as in \cite[\S4.8]{Kac}. 
Then $V(\varpi_1)$ is the natural representation of
$\bu_q(\lie g)$ and for the rest of the section we set $\varpi=\varpi_1$. 
Moreover, it is known (cf. say~\cite{CKirres}) that
if we define an $\ell$-tuple of polynomials $\bvpi=(\pi_i(u))_{i\in I}$ by
$\pi_i(u)=\delta_{i,1}(1-u)$, then $V(\bvpi)\cong V(\varpi)$ as $U_q(\lie g)$-modules.
\begin{lem} 
We have $\dim V(\bvpi)_\mu=1$ for all $\mu\in\Omega(V(\bvpi))$ and  hence $k(\bvpi)=1$. Moreover,
for all~$v\in V(\bvpi)$, $i\in\wh I$,
\begin{align*}
E_i^2 v =\, &0 = F_i^2 v 
\intertext{if $\lie g$ is of type $A_\ell,C_\ell$ or $D_\ell$, and}
E_i^{2+\delta_{i,\ell}} v = \,&0 = F_i^{2+\delta_{i,\ell}} v 
\end{align*}
if $\lie g$ is of type $B_\ell$.
\end{lem}
\begin{pf} To prove the first statement, it is enough to note that by~\cite{L} the corresponding statements hold 
for the $\bu_q(\lie g)$-module $V(\varpi)$. If $i\in I$, then the second statements also follow for the same reason. 
If $i=0$, then the result follows  from~\lemref{REP28}.
\end{pf}

Notice that 
the conditions of~\thmref{thmB} are not satisfied for the module~$L(V(\bvpi))$
and~$\Lambda\in\wh P^+$ which is not a multiple of~$\delta$.
Indeed, the first condition of~\thmref{thmB}
reads~$2(\Lambda,\delta)<(\Lambda,\alpha_1)+1$ or, equivalently,
$(\Lambda,\delta-\alpha_{1})+(\Lambda,\delta)<1$ which is a contradiction
unless~$\Lambda\in\bz\delta$. On the other hand, $w_\circ
\varpi_{1}=-\varpi_{r'}$ for some~$r'\in I$ and so the second condition
of~\thmref{thmB} yields~$(\Lambda,\delta)+(\Lambda,\alpha_{r'})<1$ which
is impossible if~$\Lambda\notin\bz\delta$.
 
Recall the function  $n: V(\bpi)\to \bn$ defined in \secref{REP}. Since the weight spaces of $V(\bvpi)$ are 
one-dimensional, it is convenient to think of this as  a function from $\Omega(V(\bvpi))\to \bn$.
We continue to denote this function by $n$.

The main result of this section is the following
\begin{thm}\label{thmE} Let  $\Lambda\in\wh P^+$ and assume that $\Lambda$ is not a multiple of $\delta$. 
Then the filtration~$\Cx_n$ on~$X(\Lambda)\tensor L(V(\bvpi))$ is strictly decreasing.
Further,
 \begin{enumerit}
\item suppose that $\lie g$ is of type $A_\ell$, $C_\ell$ or $D_{\ell}$. Then,
$$\Cx_n/\Cx_{n+1}\cong 
\bigoplus_{\mu\in\Omega_\Lambda(\varpi)} X(\Lambda+\mu+( n+n(\mu))\delta),$$
where $\Omega_\Lambda(\bvpi)=\{\mu\in\Omega(V(\bvpi))\,:\,
 \Lambda+\mu\in\wh P^+\}$.

\item  If  $\lie g$ is of type $B_\ell$,  set
$$
\Omega_\Lambda(\bvpi)=\begin{cases}
\{ \mu\in\Omega(V(\bvpi))\,:\,
 \Lambda+\mu\in\wh P^+\},&\Lambda(\alpha_\ell^\vee) > 0\\
\{\mu\in\Omega(V(\bvpi))\setminus\{0\}\,:\,
 \Lambda+\mu\in\wh P^+\},&\Lambda(\alpha_\ell^\vee)= 0.
\end{cases}
$$
Then
$$\Cx_n/\Cx_{n+1}\cong 
\bigoplus_{\mu\in\Omega_\Lambda(\varpi)} X(\Lambda+\mu+(n+n(\mu))\delta).$$
\end{enumerit}
\end{thm}
We prove this theorem in the rest of the section.
By~\propref{FIL35} it is enough to consider the case $n=0$.

\subl{RN20}
We will need the following
\begin{lem}
Suppose that there exist a sequence of integrable $\whbu_q$-modules~$\Cv_{0}\supseteq \Cv_{1}\supseteq
\cdots\supseteq \Cv_{K}$ such that $\Cv_0/\Cv_K$ is a module in the
category~$\Co$ and either~$\Cv_{i}=\Cv_{i+1}$
or~$\Cv_{i}/\Cv_{i+1}\cong X(\mu_i)$ for some~$\mu_i\in\wh P^+$. 
Set~$J=\{0\le i<K\,:\, \Cv_{i}\supsetneq \Cv_{i+1}\}$ and suppose that~$J$ is not empty.
Then $\Cv_0/\Cv_{K}\cong \bigoplus_{i\in J} X(\mu_i)$.
\end{lem}
\begin{pf}
We argue by induction on the cardinality of~$J$. Suppose first that~$J=\{i\}$
for some~$0\le i\le K$. Then
$$
\Cv_{K}=
\cdots =\Cv_{i+1}\subsetneq \Cv_{i}=\Cv_{i-1}=\cdots=\Cv_0
$$ 
and so $\Cv_0/\Cv_{K}=\Cv_{i}/\Cv_{i+1}\cong  X(\mu_i)$.

If  the cardinality of~$J$ is greater than~$1$, let~$i_1$
be the minimal  element of~$J$, that is. $\Cv_{i_1}=\Cv_{i_1-1}=\cdots =\Cv_0$ and
$\Cv_{i_1+1}\subsetneq \Cv_{i_1}$. 
 Then we have the following short exact sequence
$$
0\longrightarrow \Cv_{i_1+1}/\Cv_K\longrightarrow \Cv_0/\Cv_K
\longrightarrow \Cv_{0}/\Cv_{i_1+1}\longrightarrow 0.
$$
All modules involved are in the category~$\Co$ and integrable. Therefore, this
short exact sequence splits and so~$\Cv_0/\Cv_{K}=\Cv_{0}/\Cv_{i_1+1}
\oplus \Cv_{i_1+1}/\Cv_K
\cong  X(\mu_{i_1})\oplus\Cv_{i_1+1}/\Cv_{K}$. The lemma follows by applying
the induction hypothesis to the
 sequence~$\Cv_{K}\subseteq \cdots \subseteq\Cv_{i_1+1}$.
\end{pf}

\subl{RNA} 
Let $\lie g$ be of type $A_\ell$ and let $v_0$ be a highest weight vector in
$V(\varpi)$. Set
$$
w_0=v_0,\quad w_1=E_0 w_0,\quad w_j=E_{\ell-j+2} w_{j-1},
\qquad
2\le j\le \ell+1.
$$ 
Then~$w_{\ell+1}=w_0$ and the elements $w_j$, $0\le j\le\ell$ 
form a basis of the 
$\wh\bu_q'$-module~$V(\bvpi)$. Set~$F_{\ell+1}=F_0$. It is easy to see  
that
$$
E_jw_i=\delta_{j,\ell-i+1}w_{i+1},\qquad F_j w_{i+1}=\delta_{j,\ell-i+1} w_{i},\qquad j\in\wh I,\ \ 0\le j\le \ell, 
$$
and $n(w_j)=1-\delta_{j,0}$.
The elements $w_jt^n$, $0\le j\le\ell$, $n\in\bz$,
form a basis of the $\wh\bu_q$-module~$L(V(\bpi))$ and we have
\begin{equation}\lbl{5}
E_jw_it^n=\delta_{j,\ell-i+1}w_{i+1}t^{n+\delta_{j,0}},\qquad
F_j w_{i+1}t^n=\delta_{j,\ell-i+1} w_{i-1}t^{n-\delta_{j,0}}.
\end{equation}
Define~$\Cx_{n,j}=\whbu_q(v_\Lambda\tensor w_j t^{n+1-\delta_{j,0}})$, $0\le j\le\ell+1$. 
Then~$\Cx_{n,\ell+1}=\Cx_{n+1,0}$.
\begin{lem}
\begin{enumerit}
\item
For all $n\in\bz$, we have
$$\Cx_{n,0}\supseteq\Cx_{n,1}\supseteq\cdots\supseteq\Cx_{n,\ell}\supseteq\Cx_{n+1,0}.$$
Further, $\Cx_{n,j}\supsetneq \Cx_{n,j+1}$ implies 
$$\Cx_{n,j}/\Cx_{n,j+1}\cong X(\Lambda+\wt w_j+(n+1-\delta_{j,0})\delta).$$
\item  For all $i\in \wh I$, $0\le j\le \ell$, we have
$$F_i^{\Lambda_{i}+1}(v_\Lambda\otimes w_{j+1}t^n)=a_i\delta_{i,\ell-j+1}F_i^{\Lambda_i}(v_\Lambda\otimes  
w_jt^{n-\delta_{i,0}}),$$
for some $a_i\in\bc(q)^\times$.
\end{enumerit}
\end{lem}
\begin{pf} Part (i) is immediate from~\loceqref{5}. Part (ii) follows from~\loceqref{5} as well by 
applying~\lemref{RDC10} and~\lemref{RN10}.
\end{pf} 
Applying~\lemref{RN20} we conclude that
$$\Cx_{n}/\Cx_{n+1}\cong\bigoplus_{0\le j\le \ell\,:\,
\Cx_{n,j}\not=\Cx_{n,j+1}} \mskip-20mu X(\Lambda+\wt w_j+(n+n(w_j))\delta).
$$
Thus, in this case~\thmref{thmE} is equivalent to the following
\begin{prop}
For all~$0\le j\le\ell$, 
$\Cx_{0,j}=\Cx_{0,j+1}$ if and only if~$\Lambda_{\ell-j+1}=0$, where 
$\Lambda_{\ell+1}=\Lambda_0$. 
\end{prop}
\begin{pf}
If $\Lambda_{\ell-j+1}=0$,  then~$F_{\ell-j+1}\in\Ann_{\whbu_q^-}v_\Lambda$
by~\propref{REP10}(ii). Therefore, by (ii) of the above Lemma,
$v_\Lambda\otimes w_jt^{1-\delta_{j,0}}= c F_{\ell-j+1}(
v_\Lambda\otimes w_{j+1}t)$ for some~$c\in\bc(q)^\times$ and we are done.

For the converse, suppose that $v_\Lambda\otimes w_jt^{1-\delta_{j,0}} \in
\Cx_{0,j+1}$. It follows from~\loceqref{5} 
that we can write,
\begin{equation}\lbl{16} v_\Lambda\otimes w_jt^{1-\delta_{j,0}}
=\sum_{i=j+1}^{\ell+1}
Y_i(v_\Lambda\otimes w_i t) + 
\sum_{r=2}^R\sum_{j=0}^\ell X_{j,r}(v_\Lambda\otimes w_jt^r).
\end{equation}
for some $Y_i, X_{j,r}\in\wh\bu_q^-$. 
We first prove  that there exists an expression of the above form
in which the second  term is zero. Indeed,
suppose that the  second term in the right hand side is always non-zero
and let $R$ be  minimal such that an expression of the 
form~\loceqref{16} exists.
Let $j_0$ be  such that 
$X_{j_0,R}(v_\Lambda\otimes w_{j_0}t^{R})\ne 0$ and assume 
that $\varpi-\wt w_{j_0}$ is minimal  with 
this property.
Then by~\propref{FIL50}  $X_{j_0,R}\in\Ann_{\whbu_q^-} v_\Lambda$, say $X_{j_0,R}=\sum_{i\in\wh I}y_iF^{\Lambda_i+1}$ for some $y_i\in\whbu_q^-$.
If $j_0=1$, then it follows from~\lemref{RNA} that 
 $$X_{j_0,R}(v_\Lambda\tensor
w_{j_0}t^R)=y_0F_0^{\Lambda_0}(v_\Lambda
\tensor w_0 t^{R-1})$$ 
Since this is impossible by the choice of~$R$ we get  $j_0> 1$.
But then~\lemref{RNA} implies that  we get
 another expression of the form~\loceqref{16} 
with the minimal value of $\varpi-\wt w_j$ strictly greater than~$\varpi-\wt w_{j_0}$.
Repeating, we eventually obtain an expression of the form~\loceqref{16} 
where the minimal value of~$\varpi-
\wt w_{j}$ such that~$X_{j,R}(v_\Lambda\tensor w_j t^R)\not=0$
 is attained for~$j=1$
which is a contradiction.

Thus, we can write
\begin{equation}\lbl{17} 
v_\Lambda\otimes w_jt^{1-\delta_{j,0}}= 
\sum_{i=j+1}^{\ell+1}Y_i(v_\Lambda\otimes w_{i}t).
\end{equation} 

Let $i_0>j$ be maximal such that an expression of the above from exists and
 $Y_{i_0}(v_\Lambda\otimes w_{i_0}t)\ne 0$.
Then $\varpi-\wt v_{i_0}$ 
is minimal  with this property since 
$i_0> 0$. Hence  $Y_{i_0}\in \Ann_{\whbu_q^-} v_\Lambda$ by~\propref{FIL50}.
If $i_0=j+1$, then 
$Y_{j+1}$ is of weight $-\alpha_{\ell-j+1}$. It follows that
$Y_{j+1}=aF_{\ell-j+1}$ for some $a\in\bc(q)^\times$. Thus, $F_{\ell-j+1}\in
\Ann_{\whbu_q^-}v_\Lambda$ whence $\Lambda_{\ell-j+1}=0$. 
If $i_0>j+1$, then $Y_{i_0}$ is of weight $-(\alpha_{\ell-i_0+1}
+\cdots +\alpha_{\ell-j+1})$. 
Since~$F_j w_{i_0}=0$ unless~$j=\ell-i_0+2$, we conclude that~$Y_{i_0} (v_\Lambda\tensor
w_{i_0}t)=y F_{\ell-i_0+2}(v_\Lambda\tensor w_{i_0}t)$ and~$\Lambda_{\ell-i_0+2}=0$.
Thus, we get an expression of the form~\loceqref{17} where the maximal~$i>j$ such
that~$Y_i(v_\Lambda\tensor w_i t)\ne 0$ is at most~$i_0-1$. Repeating the argument, we get
to the case~$i_0=j+1$ which has already been  considered.
\end{pf}

\subl{RNC} Suppose that $\lie g$ is of type $C_\ell$. Then $n(\bvpi)=1$. 
Let~$v_0$ be a highest weight vector of~$V(\varpi)$ and set
\begin{alignat*}{2}
&w_0=v_0,\\
&w_j=E_{j-1} w_{j-1} 
,&\qquad & 1\le j\le \ell+1\\
&w_{\ell+j+1}=E_{\ell-j} w_{\ell+j},
&& 1\le j\le \ell-1.
\end{alignat*}
Then~$w_0,\dots,w_{2\ell-1}$ form a basis of~$V(\bvpi)$, $w_{2\ell}=w_0$ and~$n(w_j)=
1-\delta_{j,0}$. Set~$\Cx_{n,j}=\whbu_q(v_\Lambda\tensor w_j t^{n+n(w_j)})$, $0\le j\le 2\ell$. In
particular, $\Cx_{n,2\ell}=\Cx_{n+1,0}$. Then,
as in the case considered in~\ref{RNA}, \thmref{thmE} is equivalent to the following
\begin{prop}
For $0\le j\le \ell$,
$\Cx_{0,j}=\Cx_{0,j+1}$ if and only if~$\Lambda_j=0$. Similarly, for
$1\le j\le\ell-1$, 
$\Cx_{0,\ell+j}=\Cx_{0,\ell+j+1}$ if and only if~$\Lambda_{\ell-j}=0$.
\end{prop}
The proof repeats that of~\propref{RNA} with the obvious changes of notations
and we omit the details.

\subl{RNB} 
Let~$\lie g$ be of type $B_\ell$.
In this case the situation is somewhat more complicated since $n(\bvpi)=2$.
Let~$v_0$ be a highest weight vector of~$V(\bvpi)$ and set
\begin{alignat*}{2}
&w_0=v_0,\quad w_1=E_0w_0,\\  
&w_j=E_j w_{j-1},&\qquad & 2\le j\le \ell\\
&w_{\ell+j+1} = E_{\ell-j} w_{\ell+j},&& 0\le j\le \ell-2\\
&w_{2\ell}=E_0 w_{2\ell-1}
\end{alignat*}
Then~$w_0,\dots,w_{2\ell}$ form a basis of~$V(\bvpi)$. Set~$w_{2\ell+1}=w_0$.
We also have
\begin{alignat*}{2}
&F_iw_0=\delta_{i,1}a_0w_{2\ell-1},\quad F_iw_1= \delta_{i,1}a_1w_{2\ell}+\delta_{i,0}w_0,\\  
&F_iw_j=\delta_{i,j} w_{j-1},&\qquad & 2\le j\le \ell\\
&F_{i}w_{\ell+j+1} = \delta_{i,\ell-j} w_{\ell+j},&& 0\le j\le \ell-2\\
&F_iw_{2\ell}=\delta_{i,0}w_{2\ell-1}
\end{alignat*}
One can easily check that~$n(w_j)=1-\delta_{j,0}+\delta_{j,2\ell}$.
Define~$\Cx_{n,j}=\whbu_q(v_\Lambda\tensor w_j t^{n+n(w_j)})$, $j=0,\dots,2\ell$.

We have the following analog of~\lemref{RNA}.
\begin{lem}
\begin{enumerit}
\item  For all $n\in\bz$, we have
$$\Cx_{n,0}
\supseteq\Cx_{n,1}
\supseteq\cdots\supseteq\Cx_{n,2\ell-1}\supseteq \Cx_{n,2\ell}+\Cx_{n+1,0}\supseteq \Cx_{n+1,0}.
$$ 
Furthermore, $\Cx_{n,j}\supsetneq \Cx_{n,j+1}$, $0\le j\le 2\ell-2$ 
implies 
\begin{align*}\Cx_{n,j}/\Cx_{n,j+1}& \cong 
X(\Lambda+\wt w_j+(n+n(w_j))\delta).\\
\intertext{Similarly, $\Cx_{n,2\ell-1}\supsetneq \Cx_{n,2\ell}+\Cx_{n+1,1}$ implies}
\Cx_{n,2\ell-1}/(\Cx_{n,2\ell}+\Cx_{n+1,1})& \cong 
X(\Lambda+\wt w_{2\ell-1}+(n+n(w_{2\ell-1}))\delta),\\
\intertext{whilst~$\Cx_{n,2\ell}+\Cx_{n+1,0}\supsetneq \Cx_{n+1,0}$ implies}
(\Cx_{n,2\ell}+\Cx_{n+1,0})/\Cx_{n+1,1}&
\cong X(\Lambda+\wt w_{2\ell}+(n+n(w_{2\ell}))\delta).
\end{align*}

\item For all $i\in\wh I$, we have
\begin{align*} 
&F_i^{\Lambda_i+1}(v_\Lambda\otimes w_0t^n)= 
a_0 \delta_{i,1}F_i^{\Lambda_i}
(v_\Lambda\otimes w_{2\ell-1}t^{n}),\\ 
&F_i^{\Lambda_i+1}(v_\Lambda\otimes w_1t^n)= a_1 \delta_{i,1}F_i^{\Lambda_i}
(v_\Lambda\otimes w_{2\ell-1}t^{n}) +\delta_{i,0}F_i^{\Lambda_i}
(v_\Lambda\otimes w_{0}t^{n-1}),\\
&F_i^{\Lambda_i+1}(v_\Lambda\otimes w_jt^n)= 
\delta_{i,j}F_i^{\Lambda_i}
(v_\Lambda\otimes w_{j-1}t^{n}),\qquad 2\le j\le \ell\\
&F_i^{\Lambda_i+1}(v_\Lambda\otimes w_{\ell+j+1}t^n)= 
\delta_{i,\ell-j}F_i^{\Lambda_i}
(v_\Lambda\otimes w_{\ell+j}t^{n}) +
 \delta_{j,0}\delta_{i,\ell}F_i^{\Lambda_i-1}
(v_\Lambda\otimes w_{\ell-1} t^n),\\ 
&\mskip376mu 0\le j\le \ell-2\\
&F_i^{\Lambda_i+1}(v_\Lambda\otimes w_{2\ell}t^n)= 
\delta_{i,0}F_i^{\Lambda_i}
(v_\Lambda\otimes w_{2\ell-1}t^{n-1}). 
\end{align*}
 \end{enumerit}
\end{lem}
Thus, \thmref{thmE} reduces to the following
\begin{prop}
\begin{alignat*}{3}
&\Cx_{0,0}=\Cx_{0,1}&\iff&\Lambda_0=0\\
&\Cx_{0,j-1}=\Cx_{0,j}&\iff&\Lambda_j =0,&\qquad& 2\le j\le \ell,\\
&\Cx_{0,\ell+j+1}=\Cx_{0,\ell+j}&\iff&\Lambda_{\ell-j}=0,&& 0\le j \le \ell-2,\\&\Cx_{0,2\ell-1}=\Cx_{0,2\ell}+\Cx_{1,0}&\iff & \text{$\Lambda_0=0$ or $\Lambda_1 =0$},\\
&\Cx_{0,2\ell}+\Cx_{1,0}=\Cx_{1,0}&\iff & \Lambda_1 = 0.
\end{alignat*}
\end{prop}
\begin{pf}
The only if direction follows in all cases from~\propref{REP10}(ii) and the formulae
in (ii) of the above Lemma.

For the converse, we consider three separate cases.

\

{\bfseries Case 1.} $\Cx_{0,j}=\Cx_{0,j+1}$, $0\le j<2\ell -1$.

\ 

We can write
$$
v_\Lambda\otimes w_jt^{1-\delta_{j,0}}=\sum_{i=j+1}^{2\ell-1}Y_i(v_\Lambda\otimes w_i)+
Y_{2\ell+1}(v_\Lambda\otimes w_{2\ell+1}t)+\sum_{r=2}^R \sum_{i=0}^{2\ell}X_{i,r}(v_\Lambda\otimes w_it^r).
$$  
Arguing exactly as in~\propref{RNA}, but using~\lemref{RNB} instead, we conclude
that  there exists an expression of the form
\begin{equation}\lbl{17}v_\lambda\otimes w_jt^{1-\delta_{j,0}}=
\sum_{i=j+1}^{2\ell-1}Y_i(v_\Lambda\otimes w_it) +Y_{2\ell+1}
(v_\Lambda\otimes w_{2\ell+1}t).
\end{equation}
Let $i_0>j$ be maximal such that an expression of the above from exists and
 $Y_{i_0}(v_\Lambda\otimes w_{i_0}t)\ne 0$.
Then $\varpi-\wt v_{i_0}$ 
is minimal  with this property since 
$i_0> 0$. Hence  $Y_{i_0}\in \Ann_{\whbu_q^-} v_\Lambda$ by~\propref{FIL50}.
If $i_0=j+1$, then 
$Y_{j+1}$ is of weight $-\alpha_{j+1}$ if $2\le j<\ell$ and 
of weight $-\alpha_{2\ell-j+1}$ if $j\ge 1$.  It follows that
$Y_{j+1}=aF_{j+1}$ (respectively, $aF_{2\ell-j+1}$) for some $a\in\bc(q)^\times$ and hence 
$\Lambda_{j+1}=0$ (respectively, $\Lambda_{2\ell-j+1}=0$). 
Suppose that~$i_0 > j+1$ and set
$$
k = \begin{cases}
i_0,& 2\le i_0 \le \ell\\
2\ell-i_0+1,& \ell < i_0 \le 2\ell-1\\
1,& i_0=2\ell+1
\end{cases}
$$
Then~$F_i w_{i_0}=0$ unless~$i=k$. Therefore, $Y_{i_0}(v_\Lambda \tensor w_{i_0}t)=
y F_k^{\Lambda_{k}+1}(v_\Lambda\tensor w_{i_0}t)$ for some~$y\in\whbu_q^-$. Using~\lemref{RDC10} and the
formulae in~\lemref{RNB}(ii) we obtain an expression
of the form~\loceqref{17} where the maximal~$i>j$ such that~$Y_i(v_\Lambda\tensor w_it)
\ne 0$ is at most~$i_0-1$ if~$i_0\le 2\ell-1$ and at most~$2\ell-1$ if~$i_0=2\ell+1$.
Repeating this argument we reduce to the case~$i_0=j+1$ which has already been considered.

\ 

{\bf Case 2.} $\Cx_{0,2\ell-1}=\Cx_{0,2\ell}+\Cx_{1,0}$

\

In this case we should prove 
that either  $\Lambda_0=0$ or $\Lambda_1=0$.
Suppose that there exists an expression,
$$
v_\Lambda\otimes w_{2\ell-1}t= Y(v_\Lambda\otimes w_0t)+
\sum_{r=2}^R\sum_{i=0}^{2\ell}X_{i,r}(v_\Lambda\otimes w_it^r).
$$
Using~\lemref{RNB} we see  that as usual there must exist an expression of the form,
$$v_\Lambda\otimes w_{2\ell-1}t= Y_1(v_\Lambda\otimes w_0t)+Y_2(v_\Lambda\otimes w_{2\ell}t^2),
$$ for some $Y_1,Y_2\in\whbu_q^-$.
If $Y_2(v_\Lambda\otimes w_{2\ell}t^2)\ne 0$, then by~\propref{FIL50} we get that $Y_2\in\Ann_{\whbu_q^-} v_\Lambda$. On the other hand since $\wt w_{2\ell-1}t-\wt w_{2\ell}t^2=-\alpha_0$, we see that $Y_2=aF_0\in\Ann_{\whbu_q^-}v_{\Lambda}$ for some $a\in\bc(q)^\times$. Hence $\Lambda_0=0$ and we are done. 
Otherwise $Y_2(v_\Lambda\otimes w_{2\ell}t^2)=0$ and then $Y_1\in\Ann_{\wh\bu_q^-}v_\Lambda$. Again since $Y_1$ has weight $-\alpha_1$, it follows that $\Lambda_1=0$ and the proof of case 2 is complete.

\ 

{\bf Case 3.} $\Cx_{0,2\ell}+\Cx_{1,0}=\Cx_{1,0}$.

\

In this case we can write
$$
v_\Lambda\tensor w_{2\ell}t^2 = \sum_{i=0}^{2\ell-1} Y_i (v_\Lambda\tensor w_i t^2) +
Y_{2\ell+1}(v_\Lambda\tensor w_{2\ell+1}t^2)+
\sum_{r=3}^R \sum_{i=0}^{2\ell} X_{i,r} (v_\Lambda\tensor w_i t^r),
$$
for some~$Y_i, X_{i,r}\in\whbu_q^-$. Observe first that~$Y_0\in\whbu_q^-$ must
be of weight~$\wt w_{2\ell}t^2-\wt w_0 t\in\alpha_0-Q^+\notin-\wh Q^+$. Thus, $Y_0=0$.
Furthermore, using~\lemref{RNB}, we can reduce
the above expression to
\begin{equation}\lbl{27}
v_\Lambda\tensor w_{2\ell}t^2 = \sum_{i=1}^{2\ell-1} Y_i (v_\Lambda\tensor w_i t^2)+
Y_{2\ell+1}(v_\Lambda\tensor w_{2\ell+1}t^2).
\end{equation}
Let~$i_0$ be maximal~$i$ such that~$Y_i(v_\Lambda\tensor w_i t^2)\not=0$. Then~$\varpi-\wt v_{i_0}<
\varpi_i-\wt v_i$ for all~$i<i_0$, $i\not=2\ell$ and so~$Y_{i_0}\in\Ann_{\whbu_q^-}v_\Lambda$
by~\propref{FIL50}. Suppose first that~$i_0=1$. Then~$Y_{i_0}$ is of weight~$-\alpha_1$ and 
so~$Y_{i_0}=a F_1$ for some~$a\in\bc(q)^\times$. Then~$\Lambda_1=0$ by~\propref{REP10}(ii).
To complete the proof, it remains to observe that
the case~$i_0>1$ can be reduced to the case~$i_0=1$ by an argument similar to the one
in Case 1.
\end{pf}

\subl{RND} 
Finally, let~$\lie g$ be of type $D_{\ell}$, $\ell\ge 4$. In this case we also have $n(\bvpi)=2$.
Define
\begin{alignat*}{2}
&w_0=v_0,\quad w_1=E_0w_0,\\ 
&w_j=E_{j} w_{j-1},&\qquad & 2\le j\le \ell-1\\
&w_\ell= E_\ell w_{\ell-2},\\
&w_{\ell+1}=E_{\ell} w_{\ell-1}=E_{\ell-1} w_{\ell},\\
&w_{\ell+j}=E_{\ell-j} w_{\ell+j-1},&& 2\le j\le \ell-2\\
&w_{2\ell-1}=E_0 w_{2\ell-2}.
\end{alignat*}
Then~$w_0,\dots,w_{2\ell-1}$ form a basis of~$V(\bvpi)$. 
One can easily check that~$n(w_j)=1-\delta_{j,0}+\delta_{j,2\ell-1}$.
Define~$\Cx_{n,j}=\whbu_q(v_\Lambda\tensor w_j t^{n+n(w_j)})$, $j=0,\dots,2\ell-1$.
Then
\begin{multline*}
\Cx_{n,0}\supseteq \Cx_{n,1} \supseteq
\cdots \supseteq \Cx_{n,\ell-2}\supseteq \Cx_{n,\ell-1}+\Cx_{n,\ell}
\supseteq \Cx_{n,\ell+1}\supseteq\\
 \supseteq \Cx_{n,\ell+2}\supseteq\cdots \supseteq \Cx_{n,2\ell-2}\supseteq 
\Cx_{n+1,0}+\Cx_{n,2\ell-1}\supseteq \Cx_{n+1,0}
\end{multline*}
\thmref{thmE} is thus equivalent to the following
\begin{prop}
\begin{alignat*}{3}
& \Cx_{0,0}=\Cx_{0,1}&\iff & \Lambda_0 = 0\\
& \Cx_{0,j}=\Cx_{0,j+1}&\iff &\Lambda_{j+1}=0,&\qquad &1\le j\le \ell-3\\
& \Cx_{0,\ell-2}=\Cx_{0,\ell-1}+\Cx_{0,\ell}&\iff &\text{$\Lambda_{\ell-1}=0$ or~$\Lambda_\ell=0$}\\
& \Cx_{0,\ell-1} = \Cx_{0,\ell+1}&\iff & \Lambda_{\ell}=0\\
& \Cx_{0,\ell}=\Cx_{0,\ell+1}&\iff &\Lambda_{\ell-1}=0\\
& \Cx_{0,\ell+j-1}=\Cx_{0,\ell+j}&\iff & \Lambda_{\ell-j}=0,&& 1\le j\le \ell-2\\ 
& \Cx_{0,2\ell-2}=\Cx_{0,2\ell-1}+\Cx_{1,0}&\iff &\text{$\Lambda_0=0$ or~$\Lambda_1=0$}\\
& \Cx_{0,2\ell-1}+\Cx_{1,0}=\Cx_{1,0} & \iff &\Lambda_1 = 0
\end{alignat*}
\end{prop}
The proof is similar to that of~\propref{RNB} with the obvious changes in notations.

\begin{rem}
It is known (cf. for example~\cite{NS}) that the modules~$L(V(\bvpi))$ considered
in~\ref{RNA}--\ref{RND} admit crystal bases which in turn admit a simple
realization in the framework of Littelmann's path model. Let~$\wh\Cb(\bvpi)$ (respectively,
$\Cb(\Lambda)$) be
a subcrystal of Littelmann's path crystal isomorphic to a crystal basis of~$L(V(\bvpi))$.
Then the concatenation product $\Cb(\Lambda)\tensor \wh\Cb(\bvpi)$ contains a subcrystal
which is a disjoint union of indecomposable crystals isomorphic to~$\Cb(\Lambda+\wt b)$,
where~$b$ runs over the set of~$\Lambda$-dominant elements in~$\wh\Cb(\bvpi)$. For
the special cases considered above the two are actually isomorphic (this is proven
for the type~$A_\ell$ in~\cite{G}, but the argument given where remains valid for
the modules considered in~\ref{RNC} and~\ref{RND} and can be easily modified for the module
considered in~\ref{RNB}). Moreover, one can check that there is a bijection
between the set of~$\Lambda$-dominant elements of~$\wh\Cb(\bvpi)$ and the set~$\Omega_\Lambda(\bvpi)
\times\bz$.
\end{rem}
\newpage

\section*{List of notations}
\def\bqq{{\setbox0\hbox{$\whbu_q^+$}\setbox2\null\ht2\ht0\dp2\dp0\box2}}
\noindent
\begin{tabular}{p{1.8in}@{\bqq}l}
$I$&\ref{P10}\\
$\varpi_i$&\ref{P10}\\
$P$, $Q$, $Q^+$, $\Ht$&\ref{P10}\\
$\hatg$&\ref{P15}\\
$\wh I$&\ref{P15}\\
$\omega_i$&\ref{P15}\\
$\wh P$, $\wh Q$, $\wh Q^+$, $\Ht$&\ref{P15}\\
$\whbu_q$ & \ref{P20}\\
$x_{i,k}^\pm$, $h_{i,k}$&\ref{P20}\\
$K_i$&\ref{P20}\\
$C$, $D$&\ref{P20}\\
$\whbu_q^r(\ll)$, $\whbu_q^r(\gg)$, $\whbu_q^r(0)$&\ref{P20}\\
$\whbu_q^\circ$&\ref{P20}\\
$\phi_z$&\ref{P23}\\
$\whbu_q^+$, $\whbu_q^-$, $\whbu_q'$&\ref{P25}\\
$\Delta$&\ref{P30}\\
$\wh\Cu_q^+$, $\wh\Cu_q^-$&\ref{P30}\\
\end{tabular}
\hfill
\begin{tabular}{p{1.8in}@{\bqq}l}
$P_{i,k}$&\ref{P30}\\
$M_\mu$&\ref{REP0}\\
$\wt$&\ref{REP0}\\
$\Omega(M)$&\ref{REP0}\\
$M(\Lambda)$, $m_\Lambda$&\ref{REP10}\\
$X(\Lambda)$, $v_\Lambda$&\ref{REP10}\\
$V(\bpi)$, $v_\bpi^{\phantom{*}}$, $v_\bpi^*$&\ref{REP18}\\
$\lambda_{\bpi}$&\ref{REP18}\\
$w_\circ$&\ref{REP18}\\
$k(\bpi)$&\ref{REP25}\\
$n(v)$, $n(\bpi)$&\ref{REP35}\\
$m(\bpi)$ & \ref{REP40}\\
$\eta_\bpi$&\ref{REP40}\\
$V(\bpi)^{(k)}$&\ref{REP40}\\
$L(V(\bpi))$, $L^s(V(\bpi))$&\ref{REP50}\\
$\wh\eta_\bpi$&\ref{REP50}\\
$\Ht_\bpi$&\ref{SF30}\\
\end{tabular}

\bibliographystyle{amsplain}

\end{document}